\documentclass[12pt]{article}
\title{
Efficient robust nonparametric estimation in a semimartingale regression
model
\thanks{The paper is  supported by the RFBR-Grant 09-01-00172-a.}
}

\author{Konev Victor
\thanks{
Department of Applied Mathematics and Cybernetics,
 Tomsk State University,
Lenin str. 36,
 634050 Tomsk, Russia,
 e-mail: vvkonev@mail.tsu.ru }
 and
Pergamenshchikov Serguei\thanks{
 Laboratoire de Math\'ematiques Raphael Salem,
 Avenue de l'Universit\'e, BP. 12,
  Universit\'e de Rouen,
   F76801, Saint Etienne du Rouvray, Cedex France
and Department of Mathematics and Mechanics,Tomsk State University,
Lenin str. 36, 634041 Tomsk, Russia, e-mail:
Serge.Pergamenchtchikov@univ-rouen.fr } }

\usepackage{amssymb}
\usepackage{amsfonts}
\usepackage{amsmath}
\usepackage{amsthm}
\usepackage{enumerate}
\usepackage{multicol}

\newtheorem{theorem}{Theorem}[section]
\newtheorem{proposition}[theorem]{Proposition}
\newtheorem{lemma}[theorem]{Lemma}

\newtheorem{remark}{Remark}[section]
\newtheorem{corollary}[theorem]{Corollary}

\newcommand\cA{{\cal A}}
\newcommand\cC{{\cal C}}

\newcommand\cH{{\cal H}}
\newcommand\cG{{\cal G}}

\newcommand\cL{{\cal L}}
\newcommand\cB{{\cal B}}
\newcommand\cK{{\cal K}}

\newcommand\cP{{\cal P}}

\newcommand\cD{{\cal D}}
\newcommand\cQ{{\cal Q}}

\newcommand\cR{{\cal R}}
\newcommand\cT{{\cal T}}

\newcommand\cV{{\cal V}}

\newcommand\ov{\overline}

\def\bbr{{\mathbb R}}

\def\text#1{\hbox{#1}}
\def\proof{{\noindent \bf Proof. }}
\def\endproof{\mbox{\ $\qed$}}
\def\card{\mbox{card}}

\def\E{{\bf E}}
\def\P{{\bf P}}
\def\C{{\bf C}}
\def\D{{\bf D}}
\def\H{{\bf H}}
\def\M{{\bf M}}
\def\L{{\bf L}}

\newcommand{\wh}{\widehat}
\newcommand{\wt}{\widetilde}

\newcommand\Er{\mbox{Err}}

\def\Chi{{\bf 1}}
\def\d{\mathrm{d}}
\def\build #1_#2{\mathrel{\mathop{\kern 0pt #1}\limits_{#2}}}
\newcommand{\zs}[1]{{\mathchoice{#1}{#1}{\lower.25ex\hbox{$\scriptstyle#1$}}
{\lower0.25ex\hbox{$\scriptscriptstyle#1$}}}}
\numberwithin{equation}{section}

\begin{document}

\maketitle

\abstract{
The paper considers the problem of robust estimating a periodic
function in a continuous time regression model with
dependent disturbances given
by a general square integrable semimartingale  with unknown distribution.
An example of such a noise is non-gaussian Ornstein-Uhlenbeck process with
the L\'evy process subordinator,  which is used
to model the financial Black-Scholes type markets with
jumps. An adaptive model selection procedure, based on the weighted least square estimates,
  is proposed.
Under general moment conditions on the noise distribution,
 sharp
non-asymptotic oracle inequalities for the robust risks have been derived
 and the robust efficiency of the model selection procedure  has been shown.
}

\vspace*{5mm}
\noindent {\sl Keywords}: Non-asymptotic estimation; Robust risk;
 Model selection; Sharp oracle
inequality; Asymptotic efficiency.

\vspace*{5mm}
\noindent {\sl AMS 2000 Subject Classifications}:
 62G08, 62G05

\newpage

\section{Introduction}\label{sec:In}

Consider a regression model in continuous time
 \begin{equation}\label{sec:In.1}
  \d y_t=S(t)\d t+\d \xi_t\,,
  \quad 0\le t\le n\,,
 \end{equation}
where $S$ is an unknown $1$-periodic $\bbr\to\bbr$ function,
$S\in\cL_\zs{2}[0,1]$; $(\xi_\zs{t})_\zs{t\ge 0}$ is an
unobservable semimartingale noise with the values in the
Skorokhod space $\cD[0,n]$
 such that, for any
function $f$ from $\cL_\zs{2}[0,n]$, the stochastic integral
\begin{equation}\label{sec:In.2}
I_\zs{n}(f)=\int^n_\zs{0}f_\zs{s}\d \xi_\zs{s}
 \end{equation} is
well defined and has the following properties
\begin{equation}\label{sec:In.3}
\E_\zs{Q} I_\zs{n}(f)=0 \quad\mbox{and}\quad \E_\zs{Q}
I^2_\zs{n}(f)\le \sigma_\zs{Q} \int^{n}_\zs{0}\,f^{2}_\zs{s}\d s
\,.
\end{equation}
 Here $\E_\zs{Q}$ denotes the expectation with respect to the distribution
 $Q$ in $\cD[0,n]$ of the process $(\xi_\zs{t})_\zs{0\le t\le n}$,
 which is assumed to belong to some
 probability family $\cQ_\zs{n}$ specified below; $\sigma_\zs{Q}>0$ is
some positive constant depending on the distribution $Q$.

The problem is to estimate the unknown function $S$ in the model
\eqref{sec:In.1} on the basis of observations
$(y_\zs{t})_\zs{0\le t\le n}$.

The class of the disturbances $\xi$ satisfying conditions
\eqref{sec:In.3} is rather wide and comprises, in particular, the
L\'evy processes which are used in different important problems
(see \cite{Be}, for details). The models \eqref{sec:In.1}
with the L\'evy's type noise naturally arise
 (see \cite{KoPe3}) in the nonparametric functional statistics problems
(see, for example, \cite{FeVi}). Moreover, as is shown in Section~\ref{sec:Ex},
Non-Gaussian
Ornstein-Uhlenbeck-based models also enter this class.
The latter models are successfully used to model
 the Black-Scholes type financial markets
with jumps (see \cite{BaNi}, \cite{DeKl} for details and other references).

 We define the error of an estimate $\wh{S}$ (any real-valued function measurable with
respect to $\sigma\{y_\zs{t}\,,\,0\le t\le n\}$) for $S$ by its integral
quadratic risk
\begin{equation}\label{sec:In.4}
\cR_\zs{Q}(\wh{S},S):=
\E_\zs{Q,S}\,\|\wh{S}-S\|^2\,,
\end{equation}
where $\E_\zs{Q,S}$ stands for the expectation with respect to the
distribution $\P_\zs{Q,S}$ of the process \eqref{sec:In.1} with a fixed
distribution $Q$ of the noise $(\xi_\zs{t})_\zs{0\le t\le n}$
and a given function $S$; $\|\cdot\|$ is the norm in $\L_\zs{2}[0,1]$, i.e.
\begin{equation}\label{sec:In.5}
\|f\|^2:=\int^1_\zs{0}f^2(t) \d t\,.
\end{equation}
Since in our case the noise distribution $Q$ is unknown, it seems natural to measure the quality of an estimate $\wh{S}$
by the robust risk defined as
\begin{equation}\label{sec:In.6}
\cR^{*}_\zs{n}(\wh{S},S)=\sup_\zs{Q\in\cQ_\zs{n}}\,
\cR_\zs{Q}(\wh{S},S)
\end{equation}
which assumes taking supremum of the error
\eqref{sec:In.4} over the whole family of admissible distributions $\cQ_\zs{n}$.

%We  consider the estimation problem for the
%regression model \eqref{sec:In.1} for the square integrable robust risk
%  which is required to measure the quality of
%an estimate $\wh{S}$  provided that about  a true
%distribution of the noise $(\xi_\zs{t})_\zs{0\le t\le n}$ is known
% only that it  belongs to some family of distributions $\cQ_\zs{n}$.

 %It will be observed that the notion
%"nonparametric robust risk" was primarily introduced
% in \cite{GaPe0} for estimating a regression curve at a fixed point.
%The greatest lower bound for such risks have been derived and
%a point estimate is found for which
%this bound is attained. The latter means that the point estimate turns out to
%be robust efficient. In \cite{Br} and  \cite{GaPe3} this approach was applied for pointwise estimation in a heteroscedastic regression model.
%In this paper, similarly to \cite{GaPe3} and \cite{KoPe4}, we define

It is natural to treat the stated problem with respect to
the quadratic risk from the standpoint of the model selection approach.
It will be noted that the origin of this method goes back to early seventies with
the pioneering papers by Akaike \cite{Ak} and Mallows \cite{Ma} who proposed to introduce
penalizing in a log-likelihood type criterion.
 The further progress has been made by Barron,
 Birg\'e and Massart \cite{BaBiMa}, \cite{Mas},
who developed a non-asymptotic model selection method which enables one to derive
non-asymptotic oracle inequalities for nonparametric regression models
with the i.i.d. gaussian disturbances.
An oracle inequality yields the upper bound for the estimate risk via the minimal
risk corresponding to a chosen family of estimates.
 Galtchouk and Pergamenshchikov \cite{GaPe1} applied the Barron-Birg\'e-Massart
technic to the problem of estimating a nonparametric drift
function in ergodic diffusion processes.   Fourdrinier and Pergamenshchikov \cite{FoPe} extended the
 Barron-Birg\'e-Massart method to the models with
the spherically symmetric dependent observations.
They proposed a model selection procedure based on the improved least squares estimates.
Lately, the authors \cite{KoPe2} applied this method
to  the nonparametric problem of estimating
a periodic function in a model with a gaussian  colored noise in continuous
time. In all cited papers, the non-asymptotic oracle inequalities
have been derived, which
 enable one to establish the optimal convergence
rate for the minimax risks.
In addition to the optimal convergence rate, the other important problem is that of the efficiency
of adaptive estimation procedures. In order to examine the efficiency property of a
procedure  one has to obtain
 the {\sl sharp oracle inequalities}, i.e.
such in which
  the factor at the principal term in the right-hand of the inequality
 is close to unity.

The first result on sharp inequalities is most likely due to Kneip \cite{Kn}
who studied a gaussian regression model.
%The oracle
%inequalities of this type were obtained as well in \cite{CaGo} and in
%\cite{CaGoPiTs} for the inverse problems.
It will be observed that the
derivation of oracle inequalities usually rests upon
the fact that the initial model, by applying the Fourier transformation,
is reduced to the gaussian model with independent observations.
However, such a transform is possible only for
gaussian models with independent homogeneous observations or for
the inhomogeneous ones with the known correlation characteristics.
This restriction significantly narrows the area of application of the proposed
model selection procedures and rules out a broad class of models
including, in particular, widely used in econometrics
heteroscedastic regression models (see, for example, \cite{GoQu}).
For constructing adaptive procedures in the case of inhomogeneous
observations one needs to modify the approach to the estimation
problem. Galtchouk and Pergamenshchikov
 \cite{GaPe2}-\cite{GaPe3} have developed a new estimation method intended for the heteroscedastic
 regression models in discrete time. The heart of this method is to combine the
 Barron-Birg\'e-Massart non-asymptotic penalization method \cite{BaBiMa} and the
 Pinsker weighted least square method
minimizing the asymptotic risk
  (see, for example, \cite{Nu}, \cite{Pi}).
  This yields a significant improvement
in the performance of the procedure (see numerical example in \cite{GaPe2}).

The goal of this paper is to develop the robust efficient model selection method
for the model \eqref{sec:In.1} with dependent disturbances having unknown distribution.
We follow the approach, proposed by Galtchouk and Pergamenshchikov in \cite{GaPe2},
 in the construction of the procedure. Unfortunately, their method of obtaining the oracle
 inequalities is essentially based on the independence of observations and
 can not be applied here. The paper proposes the new analytical tools
 which allow one to obtain the sharp non-asymptotic oracle inequalities
for robust risks  under general conditions on the distribution of the noise
 in the model \eqref{sec:In.1}. This method enables us to treat
 both the cases of dependent and independent observations from the same standpoint, does not
assume the knowledge of the noise distribution and leads to the efficient estimation
procedure with respect to the risk \eqref{sec:In.6}. The validity of the conditions imposed
on the noise in the equation \eqref{sec:In.1}
 is verified for a non-gaussian Ornstein-Uhlenbeck process (see Section~\ref{sec:Ex}).

%  As was shown in
%\cite{GaPe3} and \cite{GaPe4}, the Galthouk-Pergamenshchikov
%procedure is efficient with respect to the robust risk \eqref{sec:In.6}.

%In the sequel in
%\cite{GaPe5}, \cite{GaPe6}, this approach has been applied to the
%problem of a drift estimation in a  diffusion process. In this
%paper we apply this procedure to the estimation of a regression
%function $S$ in a semimartingale regression model \eqref{sec:In.1}.

 The rest
of the paper is organized as follows. In Section~\ref{sec:Mo} we
construct the model selection procedure on the basis of
weighted least squares estimates
 and state the main results in the form of oracle inequalities for the quadratic
 risk \eqref{sec:In.4} and the robust risk \eqref{sec:In.6}. Here we specify also the set
 of admissible weight sequences in the model selection procedure.
  In Section~\ref{sec:Ou} we proof some properties of the stochastic integrals with
  respect to the non-gaussian Ornstein-Uhlenbeck process \eqref{sec:Ex.1}.
  Section~\ref{sec:Pr}  gives the proofs of the main results. In Sections~\ref{sec:Ef},
  \ref{sec:Up} it is shown that the proposed model selection procedure for estimating
  $S$ in \eqref{sec:In.1} is asymptotically efficient with respect to the robust risk
  \eqref{sec:In.6}.    In Appendix some auxiliary propositions are given.

\section{Non-Gaussian Ornstein-Uhlenbeck process}\label{sec:Ex}

In this section we consider an
important example of the disturbances $(\xi_\zs{t})_\zs{t\ge 0}$
in \eqref{sec:In.1} given by a non-gaussian Ornstein-Uhlenbeck process with
the L\'evy subordinator. Such processes are used in the financial
Black-Scholes type markets with  jumps
(see, for example \cite{DeKl} and the references therein).
Let the noise process in \eqref{sec:In.1} obey the equation
\begin{equation}\label{sec:Ex.1}
\d \xi_\zs{t} = a \xi_\zs{t}\d t+\d u_\zs{t}\,,\quad
\xi_\zs{0}=0\,,
\end{equation}
where $a\le 0$, $u_\zs{t} =\varrho_\zs{1}w_\zs{t}+\varrho_\zs{2}z_\zs{t}$,
$\varrho_\zs{1}$ and $\varrho_\zs{2}$ are unknown constants,
$(w_\zs{t})_\zs{t\ge 0}$ is
a standard Brownian motion, $(z_\zs{t})_\zs{t\ge 0}$ is a compound
Poisson process  defined as
$$
z_\zs{t}=\sum^{N_\zs{t}}_\zs{j=1}\,Y_\zs{j}\,.
$$
Here $(N_\zs{t})_\zs{t\ge 0}$ is a standard homogeneous Poisson process with
unknown  intensity $\lambda>0$
 and
$(Y_\zs{j})_\zs{j\ge 1}$ is an i.i.d. sequence of random variables with
\begin{equation}\label{sec:Ex.2}
\E Y_\zs{j}=0\,,\quad \E Y^2_\zs{j}=1
\quad\mbox{and}\quad
\E Y^{4}_\zs{j}<\infty\,.
\end{equation}
Let $(T)_\zs{k\ge 1}$ denote the arrival times of the process $(N_\zs{t})_\zs{t\ge 0}$,
that is,
\begin{equation}\label{sec:Ex.2-1}
T_\zs{k}=\inf\{t\ge 0\,:\,N_\zs{t}=k\}\,.
\end{equation}
 We assume that the parameters $\lambda$, $a$, $\varrho_\zs{1}$
and $\varrho_\zs{2}$ satisfy the conditions
\begin{equation}\label{sec:Ex.3}
-a_\zs{max}\le a\le 0\,,\quad
0\le \lambda\le \lambda_\zs{max}
\,,\quad \varrho^{*}_\zs{min}\le \varrho^{*}\le \varrho^{*}_\zs{max}\,,
\end{equation}
 where $\varrho^{*}=\varrho^{2}_\zs{1}
+\lambda \varrho^{2}_\zs{2}$.
Let $\cQ_\zs{n}$ denote
the family of all distributions of process
\eqref{sec:Ex.1} on $\cD[0,n]$ with the parameters $a$, $\lambda$,
$\varrho_\zs{1}$ and $\varrho_\zs{2}$ satisfying the conditions
\eqref{sec:Ex.3} with fixed  bounds $\lambda_\zs{max}>0$,
$a_\zs{max}>0$, $\varrho^{*}_\zs{min}>0$ and $\varrho^{*}_\zs{max}>0$.

\section{Model selection}\label{sec:Mo}

This Section gives the construction of a model selection procedure  for
estimating a function $S$ in \eqref{sec:In.1} on the basis of weighted least square estimates and states
the main results.

For estimating the unknown function $S$ in the model \eqref{sec:In.1}, we
apply its Fourier expansion in the trigonometric basis
  $(\phi_j)_\zs{j\ge 1}$ in $\cL_2[0,1]$
defined as
\begin{equation}\label{sec:Mo.1}
\phi_1=1\,,\quad
\phi_\zs{j}(x)=\sqrt{2}\,Tr_\zs{j}(2\pi [j/2]x)\,,\ j\ge 2\,,
\end{equation}
where the function $Tr_\zs{j}(x)=\cos(x)$ for even $j$ and
$Tr_\zs{j}(x)=\sin(x)$ for odd $j$; $[x]$ denotes the integer part
of $x$. The corresponding Fourier coefficients
\begin{equation}\label{sec:Mo.2}
\theta_\zs{j}=(S,\phi_j)= \int^1_\zs{0}\,S(t)\,\phi_\zs{j}(t)\,\d
t
\end{equation}
can be estimated as
\begin{equation}\label{sec:Mo.3}
\wh{\theta}_\zs{j,n}= \frac{1}{n}\int^n_\zs{0}\,\phi_j(t)\,\d
y_\zs{t}\,.
\end{equation}
In view of \eqref{sec:In.1}, we obtain
\begin{equation}\label{sec:Mo.4}
\wh{\theta}_\zs{j,n}=\theta_\zs{j}+\frac{1}{\sqrt{n}}\xi_\zs{j,n}\,,
\quad
 \xi_\zs{j,n}=\frac{1}{\sqrt{n}}
I_\zs{n}(\phi_\zs{j})
\end{equation}
where $I_\zs{n}(\phi_\zs{j})$ is given in \eqref{sec:In.2}.

\noindent For any sequence $x=(x_\zs{j})_\zs{j\ge 1}$, we set
\begin{equation}\label{sec:Mo.5}
|x|^2=\sum^\infty_\zs{j=1}x^2_\zs{j}
\quad\mbox{and}\quad
\#(x)=\sum^\infty_\zs{j=1}\,\Chi_\zs{\{|x_\zs{j}|>0\}}\,.
\end{equation}
Now we impose some additional conditions on the
distribution of the
 noise $(\xi_\zs{t})_\zs{t\ge 0}$ in \eqref{sec:In.1}.\\[2mm]

\noindent $\C_\zs{1})$ {\it There exists a positive constant
$\varsigma_\zs{Q}> 0$ such that
for any $n\ge 1$
$$
\L_\zs{1,n}(Q)=
 \sup_\zs{x\in \cH\,,\,\#(x)\le n}\,
 \left|
\sum^\infty_\zs{j=1}\,x_\zs{j}\,
\left(
\E_\zs{Q}\,\xi^{2}_\zs{j,n}
-
\varsigma_\zs{Q}
\right)
\right|
<\infty\,,
$$
where $\cH=[-1,1]^{\infty}$.
}

\vspace*{4mm}

\noindent $\C_\zs{2})$ {\it Assume that for all $n\ge 1$
$$
\L_\zs{2,n}(Q)=
 \sup_\zs{|x|\le 1\,,\,\#(x)\le n}\,
\E_\zs{Q}\, \left(
\sum^\infty_\zs{j=1}\,x_\zs{j}\,
(
\xi^2_\zs{j,n}-\E_\zs{Q} \xi^2_\zs{j,n}
)
%\wt{\xi}_\zs{j,n}\,
\right)^2
<\infty\,.
$$
}
\vspace{2mm}

As is shown in the proof of Theorem~\ref{Th.sec:2.2} in Section~\ref{sec:Pr} ,
 both Conditions $\C_\zs{1})$
and $\C_\zs{2})$ hold for the process \eqref{sec:Ex.1}.
 Further we define a class of weighted least squares estimates for $S(t)$ as
\begin{equation}\label{sec:Mo.6}
\wh{S}_\zs{\gamma}=\sum^{\infty}_\zs{j=1}\gamma(j)\wh{\theta}_\zs{j,n}\phi_\zs{j}\,,
\end{equation}
where $\gamma=(\gamma(j))_\zs{j\ge 1}$ is
a sequence of weight coefficients such that
\begin{equation}\label{sec:Mo.6-1}
0\le \gamma(j)\le 1 \quad\mbox{and}\quad 0<\#(\gamma)\le n\,.
\end{equation}
 Let $\Gamma$ denote a finite set of weight sequences $\gamma=(\gamma(j))_\zs{j\ge 1}$ with
 these properties, $\nu=\card(\Gamma)$ be its cardinal number and
\begin{equation}\label{sec:Mo.7}
\mu=\max_\zs{\gamma\in\Gamma}\,\#(\gamma)\,.
\end{equation}
\noindent The model selection procedure for the unknown function
$S$ in \eqref{sec:In.1} will be constructed on the basis of
a family of  estimates $(\wh{S}_\zs{\gamma})_\zs{\gamma\in\Gamma}$.
The choice of a specific set
of weight sequences
 $\Gamma$ is discussed at the end
of this section. In order to find a proper weight sequence
$\gamma$ in the set $\Gamma$ one needs  to specify a cost function.
When choosing an appropriate cost function one can use the
following argument. The empirical squared error
$$
\Er_\zs{n}(\gamma)=\|\wh{S}_\zs{\gamma}-S\|^2
$$
can be written as
\begin{equation}\label{sec:Mo.8}
\Er_\zs{n}(\gamma)\,=\,
\sum^\infty_\zs{j=1}\,\gamma^2(j)\wh{\theta}^2_\zs{j,n}\,-
2\,\sum^\infty_\zs{j=1}\,\gamma(j)\wh{\theta}_\zs{j,n}\,\theta_\zs{j}\,+\,
\sum^\infty_\zs{j=1}\theta^2_\zs{j}\,.
\end{equation}
Since the Fourier coefficients $(\theta_\zs{j})_\zs{j\ge 1}$ are
unknown, the weight coefficients $(\gamma_\zs{j})_\zs{j\ge 1}$ can not be
determined by minimizing this quantity. To circumvent this
difficulty one needs to replace  the terms
$\wh{\theta}_\zs{j,n}\,\theta_\zs{j}$ by some their estimators
$\wt{\theta}_\zs{j,n}$. We set
\begin{equation}\label{sec:Mo.9}
\wt{\theta}_\zs{j,n}=
\wh{\theta}^2_\zs{j,n}-\frac{\wh{\sigma}_\zs{n}}{n}
\end{equation}
where $\wh{\sigma}_\zs{n}$ is some estimator for the quantity $\varsigma_\zs{Q}$ in
the condition $\C_\zs{1})$.

For this change in the empirical squared error, one has to pay
some penalty. Thus, one comes to the cost function of the form
\begin{equation}\label{sec:Mo.10}
J_\zs{n}(\gamma)\,=\,\sum^\infty_\zs{j=1}\,\gamma^2(j)\wh{\theta}^2_\zs{j,n}\,-
2\,\sum^\infty_\zs{j=1}\,\gamma(j)\,\wt{\theta}_\zs{j,n}\,
+\,\rho\,\wh{P}_\zs{n}(\gamma)
\end{equation}
where $\rho$ is some positive constant,
$\wh{P}(\gamma)$ is the penalty term defined as
\begin{equation}\label{sec:Mo.11}
\wh{P}_\zs{n}(\gamma)=\frac{\wh{\sigma}_\zs{n}\,|\gamma|^2}{n}
\,.
\end{equation}

\noindent In the case, when the value of $\sigma$ in $\C_\zs{1})$ is known, one can take
$\wh{\sigma}_\zs{n}=\varsigma_\zs{Q}$ and
\begin{equation}\label{sec:Mo.12}
P_\zs{n}(\gamma)=\frac{\varsigma_\zs{Q}\,|\gamma|^2}{n}\,.
\end{equation}
Substituting the weight coefficients, minimizing the cost function
\begin{equation}\label{sec:Mo.13}
\wh{\gamma}=\mbox{argmin}_\zs{\gamma\in\Gamma}\,J_n(\gamma)\,,
\end{equation}
in \eqref{sec:Mo.6} leads to the model selection procedure
\begin{equation}\label{sec:Mo.14}
\wh{S}_\zs{*}=\wh{S}_\zs{\wh{\gamma}}\,.
\end{equation}
It will be noted that
$\wh{\gamma}$ exists, since
 $\Gamma$ is a finite set. If the
minimizing sequence in \eqref{sec:Mo.13} $\wh{\gamma}$ is not
unique, one can take any minimizer.

\begin{theorem}\label{Th.sec:2.1}
Let $\cQ_\zs{n}$ be a  family of the distributions
$Q$ on $\cD[0,n]$ such that the conditions
$\C_\zs{1})$ and $\C_\zs{2})$ hold for each $Q$
from $\cQ_\zs{n}$.
 Then for any $n\ge 1$  and $0<\rho<1/3$, the estimator
\eqref{sec:Mo.14}, for each $Q\in\cQ_\zs{n}$, satisfies the oracle inequality
\begin{align}\label{sec:Mo.15}
\cR_\zs{Q}(\wh{S}_\zs{*},S)\,\le\, \frac{1+3\rho-2\rho^2}{1-3\rho}
\min_\zs{\gamma\in\Gamma} \cR_\zs{Q}(\wh{S}_\zs{\gamma},S)
+\frac{1}{n}\,\cB_\zs{Q}(n,\rho)
\end{align}
where the risk $\cR_\zs{Q}(\cdot,S)$ is defined in \eqref{sec:In.4},
$$
\cB_\zs{Q}(n,\rho)=\Psi_\zs{Q}(n,\rho)+
\frac{6\mu\,
\E_\zs{Q,S}|\wh{\sigma}_\zs{n}-\varsigma_\zs{Q}|}{1-3\rho}
$$
and
\begin{equation}\label{sec:Mo.15-1}
\Psi_\zs{Q}(n,\rho)= \frac{ 2\varsigma_\zs{Q}\sigma_\zs{Q}\nu +
4\varsigma_\zs{Q} \L_\zs{1,n}(Q) + 2\nu \L_\zs{2,n}(Q)}{\varsigma_\zs{Q}\rho(1-3\rho)}
\,.
\end{equation}
\end{theorem}
\noindent This theorem is proved in \cite{KoPe3}.\\[2mm]

Now we will obtain the oracle inequality for the model
\eqref{sec:In.1}, \eqref{sec:Ex.1}. To write down the oracle inequality in this
case, one needs the following parameters
\begin{equation}\label{sec:Mo.15-1-1}
\lambda_\zs{1}=\lambda\varrho^{2}_\zs{1}+(\lambda\varrho_\zs{2})^{2}
\quad\mbox{and}\quad
\lambda_\zs{2}=\varrho^{2}_\zs{1}\varrho^{*}+\lambda\varrho^{2}_\zs{2}\,.
\end{equation}
Moreover, we set
\begin{equation}\label{sec:Mo.15-1-2}
\M^{*}=4\varrho^{2}_\zs{1}+\varrho^{2}_\zs{2}\D^{*}_\zs{1}
+80\lambda_\zs{2}+12\D^{*}_\zs{2}+21\varrho_\zs{3}
\end{equation}
where $\D^{*}_\zs{1}=4\lambda\varrho^{2}_\zs{1}+7\lambda^{2}\varrho^{2}_\zs{2}$,
$\D^{*}_\zs{2}=4\varrho^{2}_\zs{1}\varrho^{*}+
\varrho^{2}_\zs{2}\D^{*}_\zs{1}+23\lambda_\zs{2}$
and $\varrho_\zs{3}=\lambda \varrho^{4}_\zs{2} \E Y^{4}_\zs{1}$.

\medskip

\begin{theorem}\label{Th.sec:2.2}
Let $\cQ_\zs{n}$ be a family of the distributions
of the process \eqref{sec:Ex.1} with the parameters meeting
the conditions \eqref{sec:Ex.3}.
Then, for any $n\ge 1$  and $0<\rho<1/3$, the estimator
\eqref{sec:Mo.14} satisfies, for each $Q\in\cQ_\zs{n}$, the oracle inequality
\eqref{sec:Mo.15} with
\begin{equation}\label{sec:Mo.15-1-3}
\Psi_\zs{Q}(n,\rho)= \frac{ 6\varrho^{*}_\zs{max}\nu +
4\varrho^{*}_\zs{max} \L^{*}_\zs{1} + 56\nu \M^{*}}{\varrho^{*}_\zs{min}\rho(1-3\rho)}
\,,
\end{equation}
where
$$
\L^{*}_\zs{1}=
2(1+a_\zs{max}(a_\zs{max}+1))\varrho^{*}_\zs{max}
\,.
$$
\end{theorem}
\noindent Proof of this theorem is given in Section~\ref{sec:Pr}.

\begin{remark}\label{Re.sec:Mo.1}
Note that the term \eqref{sec:Mo.15-1-3}
does not depend on the parameters specifying the distribution family $\cQ_\zs{n}$.
 This means
that the oracle inequality \eqref{sec:Mo.15} is uniform over stability region
for the process \eqref{sec:Ex.1} including the stability bound, i.e. the case
$a=0$.
\end{remark}

To obtain the oracle inequality for the robust risks
  one has to impose additional conditions
on the distribution family $\cQ_\zs{n}$ in \eqref{sec:In.6}.
To this end, we introduce the family of distributions $Q$ on
$\cD[0,n]$ with the growth restriction on $\L_\zs{1,n}(Q)+\L_\zs{2,n}(Q)$, i.e.
\begin{equation}\label{sec:Mo.15-2}
\cP^{*}_\zs{n}=\left\{Q\in\cP_\zs{n}\,:\, \L_\zs{1,n}(Q)+\L_\zs{2,n}(Q)\le
\,l_\zs{n} \right\}\,,
\end{equation}
where $\cP_\zs{n}$ denotes the family of all probability measures on
$\cD[0,n]$, $l_\zs{n}$ is a slowly increasing positive function, i.e.
$l_\zs{n}\to+\infty$ as $n\to+\infty$ and for any $\delta>0$
$$
\lim_\zs{n\to\infty}\frac{l_\zs{n}}{n^{\delta}}=0\,.
$$

\noindent In the sequel we use the following condition

\noindent $\H_\zs{0})$
{\em Assume the distribution family $\cQ_\zs{n}$ is a subset of the class
\eqref{sec:Mo.15-2}, i.e.
$\cQ_\zs{n}\subseteq \cP^{*}_\zs{n}$, such that
\begin{equation}\label{sec:Mo.15-3}
\begin{array}{rl}
0<\varsigma_\zs{*}:= \inf_\zs{Q\in\cQ_\zs{n}} \varsigma_\zs{Q} \le
\sup_\zs{Q\in\cQ_\zs{n}} \varsigma_\zs{Q}:=\varsigma^{*}&<\infty\,;
\\[4mm]
\sigma^{*}:= \inf_\zs{Q\in\cQ_\zs{n}} \sigma_\zs{Q}&<\infty\,.
\end{array}
\end{equation}
}
\vspace{2mm}
\begin{theorem}\label{Th.sec:2.3}
Assume that the family $\cQ_\zs{n}$ in the robust risk
 \eqref{sec:In.6} satisfies the condition $\H_\zs{0})$.
 Then for any $n\ge 1$  and $0<\rho<1/3$, the estimator
\eqref{sec:Mo.14}  satisfies the oracle inequality
\begin{align}\label{sec:Mo.20}
\cR^{*}(\wh{S}_\zs{*},S)\,\le\, \frac{1+3\rho-2\rho^2}{1-3\rho}
\min_\zs{\gamma\in\Gamma} \cR^{*}(\wh{S}_\zs{\gamma},S)
+\frac{1}{n}\,\cB^{*}(n,\rho)
\end{align}
where
$$
\cB^{*}(n,\rho)=\Psi^{*}(n,\rho)+
\frac{6\mu}{1-3\rho}
\sup_\zs{Q\in\cQ_\zs{n}}\,
\E_\zs{Q,S}|\wh{\sigma}_\zs{n}-\varsigma_\zs{Q}|
$$
and
$$
\Psi^{*}(n,\rho)= \frac{ 2\varsigma^{*}\sigma^{*}\nu +
4\varsigma^{*} l_\zs{n} +
2\nu l_\zs{n}}{\varsigma_\zs{*}\rho(1-3\rho)}
\,.
$$
\end{theorem}

\medskip

\subsection{Estimation of $\varsigma_\zs{Q}$}\label{subsec:Si}

Now we consider the case of unknown quantity $\sigma$ in the condition
$\C_\zs{1})$. One can estimate $\sigma$ as

\begin{equation}\label{sec:Si.1}
\wh{\sigma}_\zs{n}=\sum^n_\zs{j=l}\,\wh{\theta}^2_\zs{j,n}
\quad\mbox{with}\quad l=[\sqrt{n}]+1\,.
\end{equation}
%where $l=[\sqrt{n}]+1$.
\vspace{2mm}
\begin{proposition}\label{Pr.sec:Si.1}
Assume that the family distribution $\cQ_\zs{n}$
satisfies the condition $\H_\zs{0})$ and
the unknown function $S(\cdot)$ in the model \eqref{sec:In.1}
is continuously differentiable
for $0\le t<1$ such that
\begin{equation}\label{sec:Si.2}
 |\dot{S}|_\zs{1}=
\int^1_\zs{0}|\dot{S}(t)|\d t < +\infty\,.
\end{equation}
Then, for any $n\ge 1$,
\begin{equation}\label{sec:Si.3}
\sup_\zs{Q\in\cQ_\zs{n}}
\E_\zs{Q,S}|\wh{\sigma}_\zs{n}-\varsigma_\zs{Q}|\le
\frac{\kappa^{*}_\zs{n}(S)}{\sqrt{n}}
\end{equation}
where
$$
\kappa^{*}_\zs{n}(S)= 4|\dot{S}|^2_\zs{1}+\varsigma^{*}+
\sqrt{l_\zs{n}}+
\frac{4|\dot{S}|_\zs{1}\sqrt{\sigma^{*}}}{n^{1/4}} +
\frac{l_\zs{n}}{n^{1/2}}\,.
$$
\end{proposition}
\proof
Substituting \eqref{sec:Mo.4} in \eqref{sec:Si.1} yields
\begin{equation}\label{sec:Pr.21}
\wh{\sigma}_\zs{n}=\sum^n_\zs{j=l}\theta^2_\zs{j}
+\frac{2}{\sqrt{n}}\sum^n_\zs{j=l}\theta_\zs{j}\xi_\zs{j,n}
+\frac{1}{n}\sum^n_\zs{j=l}\xi^2_\zs{j,n}\,.
\end{equation}
Further, denoting
$$
x'_\zs{j}=\Chi_\zs{\{l\le j\le n\}}
\quad\mbox{and}\quad
x''_\zs{j}=\frac{1}{\sqrt{n}}\,\Chi_\zs{\{l\le j\le n\}}\,,
$$
we represent the last term in \eqref{sec:Pr.21} as
$$
\frac{1}{n}\sum^n_\zs{j=l}\xi^2_\zs{j,n}=
\frac{1}{n}B_\zs{1,n}(x')
+
\frac{1}{\sqrt{n}}
\,B_\zs{2,n}(x'')
+\frac{n-l+1}{n}\varsigma_\zs{Q}\,,
$$
where
$$
B_\zs{1,n}(x)=
\sum^\infty_\zs{j=1}\,x_\zs{j}\,
\left(
\E_\zs{Q,S}\,\xi^{2}_\zs{j,n}
-
\varsigma_\zs{Q}
\right)
\quad\mbox{and}\quad
B_\zs{2,n}(x)=
\sum^\infty_\zs{j=1}\,x_\zs{j}\,
(
\xi^2_\zs{j,n}-\E_\zs{Q,S} \xi^2_\zs{j,n}
)\,.
$$
Combining these equations leads to the inequality
\begin{align*}
\E_\zs{Q,S}|\wh{\sigma}_\zs{n}-\varsigma_\zs{Q}|&\le \sum_\zs{j\ge l}\theta^2_\zs{j}
+
\frac{2}{\sqrt{n}}
\E_\zs{Q,S}|\sum^n_\zs{j=l}\theta_\zs{j}\xi_\zs{j,n}|\\[2mm]
&+
\frac{1}{n}|B_\zs{1,n}(x')|
+
\frac{1}{\sqrt{n}}
\,\E_\zs{Q,S}\,|B_\zs{2,n}(x'')|
+\frac{l-1}{n}\varsigma^{*}\,.
\end{align*}
By
Lemma~\ref{Le.sec:A.3}  and
the conditions $\C_\zs{1})$, $\C_\zs{2})$, one gets
\begin{align*}
\E_\zs{Q,S}|\wh{\sigma}_\zs{n}-\varsigma_\zs{Q}|&\le \sum_\zs{j\ge l}\theta^2_\zs{j}
+
\frac{2}{\sqrt{n}}
\E_\zs{Q,S}|\sum^n_\zs{j=l}\theta_\zs{j}\xi_\zs{j,n}|\\[2mm]
&+
\frac{\L_\zs{1,n}(Q)}{n}
+
\frac{\L_\zs{2,n}(Q)}{\sqrt{n}}
+\frac{\varsigma^{*}}{\sqrt{n}}\,.
\end{align*}
In view of the inequality \eqref{sec:In.3}, the last term can be
estimated as
$$
\E_\zs{Q,S}|\sum^n_\zs{j=l}\theta_\zs{j}\xi_\zs{j,n}|
\le
\sqrt{\sigma_\zs{Q}\sum^n_\zs{j=l}\theta^2_\zs{j}}\le
\sqrt{\sigma^{*}}|\dot{S}|_\zs{1}
\frac{2}{\sqrt{l}}\,.
$$
By applying the inequalities \eqref{sec:Mo.15-3}
we obtain the upper bound \eqref{sec:Si.3}.
Hence Proposition~\ref{Pr.sec:Si.1}.

\endproof

Theorem~\ref{Th.sec:2.1} and Proposition~\ref{Pr.sec:Si.1}
imply the following result.
\begin{theorem}\label{Th.sec:Si.1}
Assume that the family distribution $\cQ_\zs{n}$
satisfies the condition $\H_\zs{0})$ and the unknown function
$S$ is continuously differentiable satisfying the condition
\eqref{sec:Si.2}.
Then, for any $n\ge 1$ and
$0<\rho<1/3$, the model selection procedure \eqref{sec:Mo.14} with the estimator
\eqref{sec:Si.1}
 satisfies the oracle
inequality
\begin{align}\label{sec:Si.4}
\cR^{*}(\wh{S}_\zs{*},S)\,\le\, \frac{1+3\rho-2\rho^2}{1-3\rho}
\min_\zs{\gamma\in\Gamma} \cR^{*}(\wh{S}_\zs{\gamma},S)
+\frac{1}{n}\,\cB^{*}_\zs{1}(n,\rho)\,,
\end{align}
where
$$
\cB^{*}_\zs{1}(n,\rho) =\Psi^{*}(n,\rho)+
\frac{6\mu\kappa^{*}_\zs{n}(S)}{(1-3\rho) \sqrt{n}}
 \,.
$$
\end{theorem}

\subsection{Specification of weights in the model selection procedure \eqref{sec:Mo.14}}\label{subsec:Ga}

We will specify the weight coefficients $(\gamma(j))_\zs{j\ge
1}$ in the way proposed in \cite{GaPe2}
 for a heteroscedastic discrete time regression
model. Consider a numerical grid of the form
\begin{equation}\label{sec:Ga.0}
\cA_\zs{n}=\{1,\ldots,k^*\}\times\{t_1,\ldots,t_m\}\,,
\end{equation}
where  $t_i=i\varepsilon$ and $m=[1/\varepsilon^2]$. We assume that both
parameters $k^*\ge 1$ and $0<\varepsilon\le 1$ are functions of $n$, i.e.
$ k^*=k^*(n)$ and $\varepsilon=\varepsilon(n)$, such that
\begin{equation}\label{sec:Ga.1}
\left\{
\begin{array}{ll}
&\lim_\zs{n\to\infty}\,k^*(n)=+\infty\,,
\quad
\lim_\zs{n\to\infty}\,\dfrac{k^*(n)}{\ln n}=0\,,\\[6mm]
&\lim_\zs{n\to\infty}\varepsilon(n)=0
\quad\mbox{and}\quad
\lim_\zs{n\to\infty}\,n^{\delta}\varepsilon(n)\,=+\infty
\end{array}
\right.
\end{equation}
for any $\delta>0$. One can take, for example,
$$
\varepsilon(n)=\frac{1}{\ln (n+1)}
\quad\mbox{and}\quad
k^*(n)=\sqrt{\ln (n+1)}\,.
$$
For each $\alpha=(\beta,t)\in\cA_\zs{n}$, we introduce the weight
sequence $\gamma_\zs{\alpha}=(\gamma_\zs{\alpha}(j))_\zs{j\ge 1}$
given as
\begin{equation}\label{sec:Ga.2}
\gamma_\zs{\alpha}(j)=\Chi_\zs{\{1\le j\le j_\zs{0}\}}+
\left(1-(j/\omega_\alpha)^\beta\right)\, \Chi_\zs{\{ j_\zs{0}<j\le
\omega_\alpha\}}
\end{equation}
where $j_\zs{0}=j_\zs{0}(\alpha)=\left[\omega_\zs{\alpha}/\ln n\right]$,
$$
\omega_\zs{\alpha}=(\tau_\zs{\beta}\,t\,n)^{1/(2\beta+1)}
\quad\mbox{and}\quad
\tau_\zs{\beta}=\frac{(\beta+1)(2\beta+1)}{\pi^{2\beta}\beta}\,.
$$
We set
\begin{equation}\label{sec:Ga.3}
\Gamma\,=\,\{\gamma_\zs{\alpha}\,,\,\alpha\in\cA_\zs{n}\}\,.
\end{equation}
It will be noted that in this case $\nu=k^* m$.

\begin{remark}\label{Re.sec:Ga.1}
It will be observed that the specific form of weights
\eqref{sec:Ga.2} was proposed by Pinsker \cite{Pi} for the
filtration problem with known smoothness of regression function
observed with an additive gaussian white noise in the continuous
time. Nussbaum \cite{Nu} used these  weights for the
gaussian regression estimation problem in discrete time.

The minimal mean square risk,  called the Pinsker
constant, is provided by the weight least squares estimate with
the weights where the index $\alpha$ depends on the smoothness
order of the function $S$. In this case the smoothness order is
unknown and, instead of one estimate, one has to use a whole family
of estimates containing in particular the optimal one.

The problem is to study the properties of the whole class of
estimates. Below we derive an oracle inequality for this class
which yields the best mean square risk up to a multiplicative and
additive constants
provided that the
 the smoothness of the unknown function $S$
is not available. Moreover, it will be shown that the
multiplicative constant tends to unity and the additive one vanishes
as $n\to\infty$ with the rate higher than any minimax rate.
\end{remark}
\noindent In view of the assumptions \eqref{sec:Ga.1}, for any $\delta>0$,
one has
$$
\lim_\zs{n\to\infty}\,\frac{\nu}{n^{\delta}}=0\,.
$$
Moreover, by \eqref{sec:Ga.2} for any $\alpha\in\cK_\zs{n}$
\begin{align*}
\sum^\infty_\zs{j=1}\,\Chi_\zs{\{\gamma_\zs{\alpha}(j)>0\}}
  \le \omega_\zs{\alpha}\,.
\end{align*}
Therefore, taking into account that $A_\zs{\beta}\le A_1<1$ for
$\beta\ge 1$, we get
$$
\mu=\mu_\zs{n}\le(n/\varepsilon )^{1/3}\,.
$$
Therefore, for any $\delta>0$,
$$
\lim_\zs{n\to\infty}\frac{\mu_\zs{n}}{n^{1/3+\delta}}=0\,.
$$
To study the asymptotic behaviour of the term
$\cB^{*}_\zs{1}(n,\rho)$ we assume that the parameter $\rho$
in the cost function  \eqref{sec:Mo.10} depends on $n$, i.e.
$\rho=\rho_\zs{n}$ such that $\rho_\zs{n}\to 0$ as $n\to\infty$
and for any $\delta>0$
\begin{equation}\label{sec:Ga.4}
\lim_\zs{n\to\infty}\, n^{\delta} \rho_\zs{n}=0\,.
\end{equation}

Applying this limiting relation to the analysis of the asymptotic
behavior of the additive term $\cD_\zs{n}(\rho)$ in
\eqref{sec:Si.4} one comes to the following result.

\begin{theorem}\label{Th.sec:La.1}
Assume that the family distribution $\cQ_\zs{n}$
satisfies the condition $\H_\zs{0})$ and the unknown function
$S$ is continuously differentiable satisfying the condition
\eqref{sec:Si.2}.
Then, for any $n\ge 1$, the model selection procedure \eqref{sec:Mo.14},
\eqref{sec:Ga.4},
\eqref{sec:Si.1}, \eqref{sec:Ga.3}
 satisfies the oracle
inequality \eqref{sec:Si.4} with the additive term $\cB^{*}_\zs{1}(n,\rho)$
obeying, for any $\delta>0$, the following limiting relation
$$
\lim_\zs{n\to\infty}\,\frac{\cB^{*}_\zs{1}(n,\rho_\zs{n})}{n^\delta}=0
 \,.
$$
\end{theorem}

\medskip

\section{Stochastic integrals with respect to the process \eqref{sec:Ex.1}}
\label{sec:Ou}

In this Section we establish some properties of a stochastic integral
\begin{equation}\label{sec:Ou.1}
I_\zs{t}(f)=\int^{t}_\zs{0}\, f_\zs{s}\d \xi_\zs{s}
\,\quad 0\le t\le n\,,
\end{equation}
with respect to to the process \eqref{sec:Ex.1}.
We will need some notations. Let us denote
\begin{equation}\label{sec:Ou.1-1}
\varepsilon_\zs{f}(t)=a\int^{t}_\zs{0}\,e^{a(t-v)}\,f(v)\,(1+e^{2av})\,\d v\,,
\end{equation}
where $f$ is $[0,+\infty)\to\bbr$ function integrated on any finite interval.
We introduce also the following transformation
\begin{equation}\label{sec:Ou.2}
\tau_\zs{f,g}(t)=
\frac{1}{2}
\int^{t}_\zs{0}
\left(2
f(s)g(s)
+
\varepsilon^{*}_\zs{f,g}(s)
\right)\,
\d s
\end{equation}
of square integrable $[0,+\infty)\to\bbr$ functions $f$ and $g$. Here
$$
\varepsilon^{*}_\zs{f,g}(t)=f(t)\varepsilon_\zs{g}(t)
+
\varepsilon_\zs{f}(t)g(t)\,.
$$
It will be noted that
\begin{equation}\label{sec:Ou.3}
a \tau_\zs{f,1}(t)=\frac{1}{2}\varepsilon_\zs{f}(t)
\quad\mbox{and}\quad
a \tau_\zs{1,1}(t)=\frac{1}{2}\left(e^{2at}-1\right)\,.
\end{equation}

\medskip

\begin{proposition}\label{Pr.sec:Ou.1}
If $f$ and $g$ are from $\cL_\zs{2}[0,n]$
then
\begin{equation}\label{sec:Ou.4}
\E\, I_\zs{t}(f)I_\zs{t}(g)=\varrho^{*}\,\tau_\zs{f,g}(t)
\end{equation}
where $\varrho^{*}$ is given
in \eqref{sec:Ex.3}.
\end{proposition}
\proof
Noting that
the process $I_\zs{t}(f)$ satisfies the stochastic equation
$$
\d I_\zs{t}(f)=a f(t)\xi_\zs{t}\d t+f(t)\d u_\zs{t}\,,\quad
I_\zs{0}(f)=0\,,
$$
and applying the Ito formula one obtains \eqref{sec:Ou.4}.
Hence Proposition~\ref{Pr.sec:Ou.1}.

\endproof

\medskip

\noindent Further, for integrated $[0,+\infty)\to\bbr$ functions $f$ and $g$,
 we define the $[0,+\infty)\times [0,+\infty)\to\bbr$ function
\begin{equation}\label{sec:Ou.5}
D_\zs{f,g}(x,z)=
\int^{x}_\zs{0}\,L^{*}_\zs{f,g}(y,z)
\,
\d y
+f(z)g(z)\,,
\end{equation}
where $L^{*}_\zs{f,g}(y,z)=g(y+z)L_\zs{f}(y,z)+f(y+z)L_\zs{g}(y,z)$;
$$
L_\zs{f}(x,z)=a e^{a x}\,\left(
f(z)+a\int^{x}_\zs{0}\,e^{a v}\,f(v+z)\,\d v
\right)\,.
$$

\medskip

\begin{proposition}\label{Pr.sec:Ou.2}
Let $\cG_\zs{k}=\sigma\{T_\zs{1},\ldots,T_\zs{k}\}$,
where  $k\ge 1$, be $\sigma$-algebra generated by the stopping times \eqref{sec:Ex.2-1},
 $f$ and $g$  be  bounded left-continuous  $[0,\infty)\times\Omega\to\bbr$ functions
 measurable with respect to $\cB[0,+\infty)\bigotimes \cG_\zs{k}$
(the product $\sigma$ algebra created by $\cB[0,+\infty)$ and $\cG_\zs{k}$).
Then
$$
\E\left( I_\zs{T_\zs{k}-}(f)
|\cG_\zs{k}
\right)
=0
$$
and
$$
\E\left( I_\zs{T_\zs{k}-}(f)\, I_\zs{T_\zs{k}-}(g)
|\cG_\zs{k}
\right)
=\varrho^2_\zs{1}
\tau_\zs{f,g}(T_\zs{k})
+
\varrho^2_\zs{2}
\sum^{k-1}_\zs{l=1}\,D_\zs{f,g}(T_\zs{k}-T_\zs{l},T_\zs{l})\,.
$$

\end{proposition}
\proof
By the Ito formula one has
\begin{align}\nonumber
I_\zs{t}(f)\, I_\zs{t}(g)&=
\int^{t}_\zs{0}\,(\varrho^{2}_\zs{1} f(s)g(s)
+a (f(s) I_\zs{s}(g) + g(s) I_\zs{s}(f))\xi_\zs{s})\d s\\[2mm] \nonumber
&+\varrho^{2}_\zs{2}\,
\sum_\zs{l\ge 1}\,f(T_\zs{l})\,g(T_\zs{l})\,Y^{2}_\zs{l}
\Chi_\zs{\{T_\zs{l}\le t\}}
\\[2mm]
\label{sec:Ou.6}&
+\int^{t}_\zs{0}\,
(f(s)I_\zs{s-}(g) + g(s)I_\zs{s-}(f)))\d u_\zs{s}\,.
\end{align}
Taking the conditional expectation $\E\left(\cdot|\cG_\zs{k}\right)$,
on the set $\{T_\zs{k}>t\}$, yields
\begin{align*}
\E \left(I_\zs{t}(f)\, I_\zs{t}(g)|\cG_\zs{k}\right)&=
\int^{t}_\zs{0}\,\varrho^{2}_\zs{1} f(s)g(s)\d s
+\varrho^{2}_\zs{2}\,
\sum_\zs{l\ge 1}\,f(T_\zs{l})\,g(T_\zs{l})\,
\,\Chi_\zs{\{T_\zs{l}\le t\}}
\\[2mm]
&+a
\int^{t}_\zs{0}\,
\left(
f(s)
\E(I_\zs{s}(g) \xi_\zs{s}|\cG_\zs{k}) + g(s)
\E(I_\zs{s}(f)\xi_\zs{s}|\cG_\zs{k})
\right)
\d s\,.
\end{align*}
Now to calculate the function $Z_\zs{t}=\E(I_\zs{t}(f)\xi_\zs{t}|\cG_\zs{k})$
we put $g=1$. Taking into account that
$$
\E ( \xi^2_\zs{t}|\cG_\zs{k})=
\frac{\varrho^{2}_\zs{1}}{2a}(e^{2at}-1)
+\varrho^{2}_\zs{2}\sum_\zs{l\ge 1}
e^{2a(t-T_\zs{l})}\,\Chi_\zs{\{T_\zs{l}\le t\}}\,,
$$
one obtains,
for $T_\zs{j-1}\le t<T_\zs{j}$,
$$
\dot{Z}_\zs{t}=aZ_\zs{t}+f(t) \psi_\zs{t}\,,
$$
where
$$
\psi_\zs{t}=\frac{\varrho^{2}_\zs{1}}{2}
\left(
1+e^{2at}
\right)+\varrho^{2}_\zs{2}a\sum_\zs{l\ge 1} e^{2a(t-T_\zs{l})}
\Chi_\zs{\{T_\zs{l}\le t\}}\,.
$$
Therefore, for $T_\zs{j-1}\le t<T_\zs{j}$,
$$
Z_\zs{t}=e^{a(t-T_\zs{j-1})} Z_\zs{T_\zs{j-1}}
+
\int^{t}_\zs{T_\zs{j-1}} e^{a(t-s)} f(s) \psi_\zs{s}\d s\,.
$$
>From here and \eqref{sec:Ou.5} with $g=1$ one has
\begin{equation}
\label{sec:Ou.7}
Z_\zs{T_\zs{j}}=e^{a(T_\zs{j}-T_\zs{j-1})} Z_\zs{T_\zs{j-1}}
+
\eta_\zs{j}\,,\quad Z_\zs{T_\zs{0}}=Z_\zs{0}=0\,,
\end{equation}
where
$$
\eta_\zs{j}=\int^{T_\zs{j}}_\zs{T_\zs{j-1}} e^{a(T_\zs{j}-s)} f(s) \psi_\zs{s}\d s
+
\varrho^{2}_\zs{2} f(T_\zs{j})\,.
$$
Solving the equation \eqref{sec:Ou.6} one obtains
$$
Z_\zs{T_\zs{j}}=\int^{T_\zs{j}}_\zs{0}\,e^{a(T_\zs{j}-s)} f(s) \psi_\zs{s}\d s
+\varrho^{2}_\zs{2}\,\sum^{j}_\zs{l=1} e^{a(T_\zs{j}-T_\zs{l})} f(T_\zs{l})\,.
$$
Therefore,
$$
Z_\zs{t}=\int^{t}_\zs{0} e^{a(t-s)} f(s) \psi_\zs{s}\d s +
\varrho^{2}_\zs{2}\sum_\zs{l\ge 1} e^{a(t-T_\zs{l})} f(T_\zs{l})
\Chi_\zs{\{T_\zs{l}\le t\}}\,,
$$
i.e
$$
a\E(I_\zs{t}(f)\xi_\zs{t}|\cG_\zs{k})=
\frac{\varrho^{2}_\zs{1}}{2}\varepsilon_\zs{f}(t)
+\varrho^{2}_\zs{2}\,\sum_\zs{j\ge 1}\,L_\zs{f}(t-T_\zs{j},T_\zs{j})\,
\Chi_\zs{\{T_\zs{j}\le t\}}.
$$
>From here one comes to the desired equality. Hence
Proposition~\ref{Pr.sec:Ou.2}.

\endproof

\medskip

\begin{proposition}\label{Pr.sec:Ou.3}
Let $F$, $f$ and $g$ be non random bounded left-continuous\\
$[0,\infty)\to\bbr$ functions.
Then
\begin{align*}
\E\,\sum_\zs{k\ge 1}\,
F(T_\zs{k})\, I_\zs{T_\zs{k-}}(f)\, I_\zs{T_\zs{k-}}(g)\,\Chi_\zs{\{T_\zs{k}\le t\}}
=\int^{t}_\zs{0}\,F(v)\,H_\zs{f,g}(v)\,\d v\,,
\end{align*}
where
$$
H_\zs{f,g}(t)=\lambda \varrho^{2}_\zs{1}\,\tau_\zs{f,g}(t)
+
(\lambda\varrho_\zs{2})^{2}\,
\int^{t}_\zs{0}\,
D_\zs{f,g}(t-z,z)
\d z
\,.
$$
\end{proposition}
\proof We have
$$
\iota(t)=\E\,\sum_\zs{k\ge 1}\,
F(T_\zs{k})\, I_\zs{T_\zs{k-}}(f)\, I_\zs{T_\zs{k-}}(g)\,\Chi_\zs{\{T_\zs{k}\le t\}}
\,.
$$
By Proposition~\ref{Pr.sec:Ou.2}
one gets
\begin{align*}
\iota(t)&=\varrho^{2}_\zs{1}\E\,\sum_\zs{k\ge 1}\,
F(T_\zs{k})\,\tau_\zs{f,g}(T_\zs{k})\,\Chi_\zs{\{T_\zs{k}\le t\}}\\[2mm]
&+
\varrho^{2}_\zs{2}\,
\E\,\sum_\zs{k\ge 1}\,
F(T_\zs{k})\,\sum^{k-1}_\zs{l=1} D_\zs{f,g}(T_\zs{k}-T_\zs{l}, T_\zs{l})\Chi_\zs{\{T_\zs{k}\le t\}}\\[2mm]
&:=\iota_\zs{1}(t)+\iota_\zs{2}(t)
\,,
\end{align*}
where
\begin{align*}
\iota_\zs{1}(t)=
\lambda
\int^{t}_\zs{0}\,
\sum_\zs{l\ge 1}\,
F(z)\tau_\zs{f,g}(z) \frac{(\lambda z)^{l-1}}{(l-1)!} e^{-\lambda z} \d z =\int^{t}_\zs{0} F(z)\tau_\zs{f,g}(z) \d z\,.
\end{align*}
To calculate $\iota_\zs{2}(t)$ we note that
\begin{align*}
\iota_\zs{2}(t)=\E \sum_\zs{l\ge 1}\,\Chi_\zs{\{T_\zs{l}\le t\}}\sum_\zs{k\ge l+1}
F(T_\zs{k})
D_\zs{f,g}(T_\zs{k}-T_\zs{l}, T_\zs{l})\Chi_\zs{\{T_\zs{k}\le t\}}
\,.
\end{align*}
Taking into account that $T_\zs{k}-T_\zs{l}$ is independent of  $T_\zs{l} $
for any $k>l$ we obtain
\begin{align*}
\iota_\zs{2}(t)&=
\lambda\E \sum_\zs{l\ge 1}\,\Chi_\zs{\{T_\zs{l}\le t\}}
\int^{t-T_\zs{l}}_\zs{0}
\sum_\zs{k\ge l+1}
F(z+T_\zs{l})
D_\zs{f,g}(z, T_\zs{l})\, \frac{(\lambda z)^{k-l-1}}{(k-l-1)!} e^{-\lambda z}
\d z\\[2mm]
&=
\lambda\E \sum_\zs{l\ge 1}\,\Chi_\zs{\{T_\zs{l}\le t\}}
\int^{t-T_\zs{l}}_\zs{0}
F(z+T_\zs{l})
D_\zs{f,g}(z, T_\zs{l})\,
\d z\\[2mm]
&=\lambda^2 \int^{t}_\zs{0}\,
\int^{t-x}_\zs{0}
\left(
F(z+x)\,D_\zs{f,g}(z,x)\d z
\right)
\d x
\,.
\end{align*}

\endproof

\medskip

\noindent Note that
\begin{equation}\label{sec:Ou.8}
a H_\zs{f,1}(t)=
\frac{\lambda_\zs{1}}{2}\,
\varepsilon_\zs{f}(t)
\quad\mbox{and}\quad
a H_\zs{1,1}(t)=\frac{\lambda_\zs{1}}{2}\left(e^{2at}-1\right)\,,
\end{equation}
where $\lambda_\zs{1}$ given in \eqref{sec:Mo.15-1-1}. Now we set
\begin{equation}\label{sec:Ou.9}
\wt{I}_\zs{t}(f)=I^{2}_\zs{t}(f)-\E\,I^{2}_\zs{t}(f)\,.
\end{equation}
\noindent Further we need the following correlation measures
for two integrated $[0,+\infty)\to \bbr$ functions $f$ and $g$
\begin{equation}\label{sec:Ou.10}
\varpi_\zs{f,g}
=\max_\zs{0\le v\le n}\,
\max_\zs{0\le t\le n-v}\,
|\int^{t}_\zs{0}f(u+v)g(u)\d u|
\end{equation}
and
\begin{equation}\label{sec:Ou.11}
\varpi^{*}_\zs{f,g}=\max\left(\varpi_\zs{f,g}, \varpi_\zs{g,f}\right)\,.
\end{equation}
For any bounded $[0,\infty)\to\bbr$ function $f$ we introduce the following uniform norm
$$
\|f\|_\zs{*}=\sup_\zs{0\le t\le n}|f(t)|\,.
$$
To check the condition $\C_\zs{2})$ one needs the following non-asymptotic upper
bound
\begin{theorem}\label{Th.sec:Ou.1}
For any bounded left-continuous $[0,\infty)\to\bbr$ functions $f$, $g$
\begin{equation}\label{sec:Ou.12}
|\E \wt{I}_\zs{n}(f)\wt{I}_\zs{n}(g)|\le n
\M^{*}\,\left(\varpi^{*}_\zs{f,g}+
\|f\|_\zs{*}\|g\|_\zs{*}\right)
\|f\|_\zs{*}\|g\|_\zs{*}
\,,
\end{equation}
where $\M^{*}$ is defined in \eqref{sec:Mo.15-1-2}.
\end{theorem}
\proof
By the Ito formula one comes to the following stochastic equation
$$
\d \wt{I}_\zs{t}(f)=
2a v_\zs{t}(f) f(t)\d t+\d M_\zs{t}(f)\,,
\quad
\wt{I}_\zs{0}(f)=0\,,
$$
with $v_\zs{t}(f)=I_\zs{t}(f)\xi_\zs{t}-
\E\,I_\zs{t}(f)\xi_\zs{t}$,
$$
\d M_\zs{t}(f)=2I_\zs{t-}(f)f(t)\d u_\zs{t}
+\varrho^{2}_\zs{2}f^{2}(t)\d m_\zs{t}\,,\quad M_\zs{0}(f)=0\,.
$$
In view of  Propositions~\ref{Pr.sec:Ou.1}--\ref{Pr.sec:Ou.3}, one finds
\begin{equation}\label{sec:Ou.13}
\E\,[\wt{I}(f),\wt{I}(g)]_\zs{t}=
\E\,[M(f),M(g)]_\zs{t}
=\int^{t}_\zs{0}\,V_\zs{f,g}(s)\,\d s\,,
\end{equation}
where
$V_\zs{f,g}(t)=4
f(t)g(t) G_\zs{f,g}(t)+ \varrho_\zs{3}\,f^2(t)\,g^{2}(t)$
and $G_\zs{f,g}=\varrho^{2}_\zs{1}\varrho^{*}
\tau_\zs{f,g}(t)
+\varrho^{2}_\zs{2}H_\zs{f,g}(t)$.
Note that Lemma~\ref{Le.sec:A.1} implies
\begin{equation}\label{sec:Ou.14}
\max_\zs{0\le t\le n}\,
|G_\zs{f,g}(t)|\,\le\,(4\varrho^{2}_\zs{1}\varrho^{*}+
\varrho^{2}_\zs{2} \D^{*}_\zs{1})\varpi^{*}_\zs{f,g}\,.
\end{equation}
One can easily check that
\begin{equation}\label{sec:Ou.15}
V_\zs{1,1}(t)=2\lambda_\zs{2}\frac{e^{2at}-1}{a}+\varrho_\zs{3}\,.
\end{equation}
The constants $\lambda_\zs{2}$,  $\varrho_\zs{3}$
and $\D^{*}_\zs{1}$ are given in \eqref{sec:Mo.15-1-1}-\eqref{sec:Mo.15-1-2}.
Moreover, by the Ito formula we get
$$
\d v_\zs{t}(f)=
a v_\zs{t}(f) \d t+ af(t)\zeta_\zs{t}\d t
+
\d K_\zs{t}(f)\,,
\quad
v_\zs{0}(f)=0\,,
$$
where $\zeta_\zs{t}=\xi^{2}_\zs{t}
-
\E\,\xi^{2}_\zs{t}$,
\begin{equation}\label{sec:Ou.16}
 K_\zs{t}(f)= \int^{t}_\zs{0}
\,I^{*}_\zs{s-}(f)\,\d u_\zs{s} +\int^{t}_\zs{0}
\varrho^{2}_\zs{2}f(s)\d m_\zs{s}
\end{equation}
with
$$
I^{*}_\zs{t}(f)
=
I_\zs{t}(f) + f(t)\xi_\zs{t}
\quad\mbox{and}\quad
m_\zs{t}=\sum_\zs{0\le s\le t}\,\Delta z^{2}_\zs{s}
-\lambda t\,.
$$
Proposition~\ref{Pr.sec:Ou.1} implies
$$
\E\,I^{*}_\zs{t}(f) I^{*}_\zs{t}(g)=\varrho^{*}\,\tau^{*}_\zs{f,g}(t)\,,
$$
where
$$
\tau^{*}_\zs{f,g}(t)=\tau_\zs{f,g}(t)+f(t)\tau_\zs{1,g}(t)
+g(t)\tau_\zs{f,1}(t)+f(t)g(t)\,\tau_\zs{1,1}(t)\,.
$$
>From \eqref{sec:Ou.2}--\eqref{sec:Ou.3} it follows
that
\begin{equation}\label{sec:Ou.17}
 \tau^{*}_\zs{f,g}(t)=\tau_\zs{f,g}(t)+
\frac{\varepsilon^{*}_\zs{f,g}(t)+f(t)g(t)(e^{2at}-1)}{2a}\,.
\end{equation}
By applying Proposition~\ref{Pr.sec:Ou.3} one finds
$$
\E\,\sum_\zs{k\ge 1}\,
 I^{*}_\zs{T_\zs{k-}}(f)\, I^{*}_\zs{T_\zs{k-}}(g)\,\Chi_\zs{\{T_\zs{k}\le t\}}
=\int^{t}_\zs{0}\,H^{*}_\zs{f,g}(v)\,\d v\,,
$$
where
$$
H^{*}_\zs{f,g}(t)=H_\zs{f,g}(t)+f(t)H_\zs{1,g}(t)
+g(t)H_\zs{f,1}(t)+f(t)g(t)\,H_\zs{1,1}(t)\,.
$$
>From here and \eqref{sec:Ou.8} we get
\begin{equation}\label{sec:Ou.18}
H^{*}_\zs{f,g}(t)=H_\zs{f,g}(t)+
\lambda_\zs{1}
\frac{\varepsilon^{*}_\zs{f,g}(t)+f(t)g(t)(e^{2at}-1)}{2a}\,.
\end{equation}
Taking this into account, we calculate that,
for any square integrated functions $f$ and $g$,
\begin{equation}\label{sec:Ou.19}
\E\,[K(f)\,,\,K(g)]_\zs{t}=
\int^{t}_\zs{0}\left(
G^{*}_\zs{f,g}(s)+\varrho_\zs{3}\,f(s)\,g(s)\,
\right)\d s\,,
\end{equation}
where
$$
G^{*}_\zs{f,g}(t)=\varrho^{2}_\zs{1}\varrho^{*}\tau^{*}_\zs{f,g}(t)
+\varrho^{2}_\zs{2}\, H^{*}_\zs{f,g}(t)
\,.
$$
Further by applying Propositions~\ref{Pr.sec:Ou.1}--\ref{Pr.sec:Ou.3} we obtain
\begin{equation}\label{sec:Ou.20}
\E\,[K(f),M(g)]_\zs{t}=\int^{t}_\zs{0}\,U_\zs{f,g}(s) \d s\,,
\end{equation}
where
$$
U_\zs{f,g}(s)=2g(s)G_\zs{f,g}(s)+2g(s)f(s)G_\zs{1,g}(s)
+\varrho_\zs{3}f(s)g^{2}(s)\,.
$$
\noindent By the Ito formula one finds for $t\ge 0$
\begin{align}\label{sec:Ou.21}
\E\,\wt{I}_\zs{t}(f)\,\wt{I}_\zs{t}(g)&=
\E\,[\wt{I}(f)\,,\,\wt{I}(g)]_\zs{t}
+2\int^{t}_\zs{0}\,
\left(
f(s)\,\cT_\zs{f,g}(s)
+
g(s)\,\cT_\zs{g,f}(s)
\right)\,\d s\,,
\end{align}
where
$\cT_\zs{f,g}(t)=a \E\,v_\zs{t}(f)\,\wt{I}_\zs{t}(g)$.
Since for $g=1$ the processes
$\wt{I}_\zs{t}(g)$ and $v_\zs{t}(g)$ coincide with $\zeta_\zs{t}$,
i.e.  $\wt{I}_\zs{t}(1)=v_\zs{t}(1)=\zeta_\zs{t}$,
  \eqref{sec:Ou.16} implies
\begin{align}\label{sec:Ou.22}
\E \zeta^{2}_\zs{t}&=\int^{t}_\zs{0}\,e^{4a(t-s)}\,V_\zs{1,1}(s)\d s
=e^{4at}
\frac{2\lambda_\zs{2}+a\varrho_\zs{3}}{4a^{2}}
+
e^{2at}
\frac{\lambda_\zs{2}}{a^{2}}
+
\frac{2\lambda_\zs{2}-a\varrho_\zs{3}}{4a^{2}}
\,.
\end{align}

\noindent We define the function
\begin{equation}\label{sec:Ou.23}
A_\zs{f}(t)=\int^{t}_\zs{0}\,e^{3a(t-s)}
\left(f(s)a^{2}\,\E\,\zeta^{2}_\zs{s}
+\kappa_\zs{f}(s)
\right)\d s\,,
\end{equation}
where
$$
\kappa_\zs{f}(t)=\lambda_\zs{2}\left(
\varepsilon_\zs{f}(t)
+f(t)\left(
e^{2at}-1
\right)
\right)
+a\varrho_\zs{3}f(t)\,.
$$
Denote
\begin{equation}\label{sec:Ou.24}
V_\zs{f,g}(s)=A^{*}_\zs{f,g}(s)
+
G^{*}_\zs{f,g}(s)+\varrho_\zs{3}\,f(s)g(s)\,,
\end{equation}
where $A^{*}_\zs{f,g}(t)=g(s) A_\zs{f}(s)
+f(s) A_\zs{g}(s)$.

\noindent To calculate the function $\cT_\zs{f,g}(t)$, we
note that, by the Ito formula and \eqref{sec:Ou.15},
\begin{align}\nonumber
 \d \E\, v_\zs{t}(f) \wt{I}_\zs{t}(g)&= a(\E v_\zs{t}(f) \wt{I}_\zs{t}(g))\d t
+2ag(t)(\E v_\zs{t}(f)\,v_\zs{t}(g))\d t\\[2mm]\label{sec:Ou.25}
&+af(t)(\E \zeta_\zs{t}\wt{I}_\zs{t}(g))\d t
+U_\zs{f,g}(t)\d t\,.
\end{align}
Substituting here $f=1$, and
then taking into account \eqref{sec:A.2} yield
$$
\E \zeta_\zs{t}\wt{I}_\zs{t}(g)=
\int^{t}_\zs{0} e^{2a(t-s)}
\left(
A^{*}_\zs{g,g}(s)
+U_\zs{1,g}(s)
\right)\d s\,.
$$
Furthermore, by \eqref{sec:A.2}
\begin{align*}
\E\,v_\zs{t}(f)\wt{I}_\zs{t}(g)&=
2a\int^{t}_\zs{0}\,e^{a(t-s)}\,g(s)\,
\left(\E v_\zs{s}(f)v_\zs{s}(g)\right)\,\d s\\[2mm]
&+a\int^{t}_\zs{0}\,e^{a(t-s)} f(s)\,
\left(\E\zeta_\zs{s}\,\wt{I}_\zs{s}(g)\right)\,\d s
%\\[2mm]&
+
\int^{t}_\zs{0}\,e^{a(t-s)}\,U_\zs{f,g}(s)\d s\,.
\end{align*}
Therefore, for any bounded left-continuous
$[0,+\infty)\to \bbr$ functions $f$ and $g$, one finds
\begin{equation}\label{sec:Ou.26}
\cT_\zs{f,g}(t)=a
\int^{t}_\zs{0}
e^{a(t-s)}
\left(
g(s)\cV_\zs{f,g}(s)
+
f(s)\cK_\zs{g}(s)
+
U_\zs{f,g}(s)
\right)
\,\d s\,,
\end{equation}
where
$$
\cV_\zs{f,g}(t)=2a\int^{t}_\zs{0}\,e^{2a(t-s)}\,V_\zs{f,g}(s)\d s
$$
and
$$
\cK_\zs{g}(t)=a
\int^{t}_\zs{0}\,e^{2a(t-s)}
\left(A^{*}_\zs{g,g}(s)+U_\zs{1,g}(s)
\right)
\d s\,.
$$
>From \eqref{sec:Ou.16} and \eqref{sec:Ou.26}, it follows that
\begin{align*}
|\E\,\wt{I}_\zs{n}(f)\,\wt{I}_\zs{n}(g)|\le
\,n
\|V_\zs{f,g}\|_\zs{*}
+2n
\|f\|_\zs{*}
\|\cT_\zs{f,g}\|_\zs{*}+
2n
\|g\|_\zs{*}
\| \cT_\zs{g,f}\|_\zs{*}\,.
\end{align*}
Now by applying the inequality \eqref{sec:A.9} one gets
$$
\|V_\zs{f,g}\|_\zs{*}
\le
\left(
(4\varrho^{2}_\zs{1}+\varrho^{2}_\zs{2}\D^{*}_\zs{1}) \varpi^{*}_\zs{f,g}
+\varrho_\zs{3} \|f\|_\zs{*} \|g\|_\zs{*}
\right)
\|f\|_\zs{*} \|g\|_\zs{*} \,.
$$
Note that Lemmas~\ref{Le.sec:A.6}--\ref{Le.sec:A.8}
imply
$$
\|\cT_\zs{f,g}\|_\zs{*}
\le \left(
(20\lambda_\zs{2}+3\D^{*}_\zs{2}) \varpi^{*}_\zs{f,g}
+5\varrho_\zs{3} \|f\|_\zs{*} \|g\|_\zs{*}
\right)\, \|g\|_\zs{*}\,.
$$
>From here one comes to the upper bound \eqref{sec:Ou.11}.

\endproof

\section{Proof of Theorem~\ref{Th.sec:2.2}}\label{sec:Pr}

First we note that Proposition~\ref{Pr.sec:Ou.1} implies
 the inequality \eqref{sec:In.3} with
$\sigma_\zs{Q}=3\varrho^{*}$. Therefore due to the conditions \eqref{sec:Ex.3}
one obtains $\sigma^{*}=3 \varrho_\zs{max}$.
 Now we  verify Conditions $\C_\zs{1})$ and $\C_\zs{2})$ for the family
of processes \eqref{sec:Ex.1} satisfying the conditions \eqref{sec:Ex.3}.
 To begin with we note that
$$
\E_\zs{Q,S} \xi^2_\zs{j,n}=
\varrho^{*}
\left(
1+
\frac{a}{n}\int^{n}_\zs{0}\,e^{av}\,\Upsilon_\zs{j}(v)\d v
\right)\,,
$$
where
$$
\Upsilon_\zs{j}(v)=\int^{n}_\zs{v}\,\phi_\zs{j}(t)\,\phi_\zs{j}(t-v)\,\left(
1+e^{2a(t-v)}
\right)\d t\,.
$$
If $j=1$, one has
\begin{equation}\label{sec:Pr.1}
|\E_\zs{Q,S} \xi^2_\zs{1,n}
-
\varrho^{*}
 |
\le 2\varrho^{*}\,.
\end{equation}
Since for the trigonometric basis \eqref{sec:Mo.1}
for $j\ge 2$
$$
\phi_\zs{j}(t)\,\phi_\zs{j}(t-v)=
\cos(\gamma_\zs{j}v)+(-)^{j}\cos(\gamma_\zs{j}(2t-v))
$$
where $\gamma_\zs{j}=2\pi [j/2]$, therefore,
$$
\Upsilon_\zs{j}(v)=\cos(\gamma_\zs{j}v) F(v)+
(-1)^{j}\,\Upsilon_\zs{0,j}(v)
\,, \quad
F(v)=\int^{n-v}_\zs{0}
\left(
1+e^{2a t}
\right)\d t
$$
 and
$$
\Upsilon_\zs{0,j}(v)=
\int^{n-v}_\zs{0}
\cos(\gamma_\zs{j}(2t+v))
\left(
1+e^{2a t}
\right)\d t\,.
$$
Integrating by parts yields
$$
\Upsilon_\zs{0,j}(v)
=\frac{e^{2a(n-v)}}{2\gamma_\zs{j}}\,
\sin(v\gamma_\zs{j})+
\frac{a}{2\gamma^{2}_\zs{j}}\,\Upsilon_\zs{1,j}(v)
$$
where
$$
\Upsilon_\zs{1,j}(v)=
\cos(v\gamma_\zs{j}) (e^{a(n-v)}-1)
-
a\int^{n-v}_\zs{0} e^{at} \cos((2t+v)\gamma_\zs{j})\,\d t\,.
$$
It is obvious, that $|\Upsilon_\zs{1,j}(v)|\le 2$.
Further we obtain
\begin{align*}
a\int^{n}_\zs{0}\,e^{av}\,&\Upsilon_\zs{j}(v)\d v=
a\int^{n}_\zs{0}\,e^{av}\,F(v)\,\cos(v\gamma_\zs{j})\d v
+
a(-1)^{j}
\int^{n}_\zs{0}\,e^{av}\,\Upsilon_\zs{0,j}(v)\,\d v
\\[2mm]
&
:=a D_\zs{1}(n)+
a(-1)^{j} D_\zs{2}(n)\,.
\end{align*}
Integrating by parts two times we find
$$
D_\zs{1}(n)=\frac{1}{\gamma^{2}_\zs{j}}
\left(
e^{an}\dot{F}(n)
-
\dot{F}(0)-aF(0)-
\int^{n}_\zs{0}\,e^{av} F_\zs{1}(v)\d v
\right)\,,
$$
where
$$
F_\zs{1}(v)=a^{2}F(v)+2a\dot{F}(v)+\ddot{F}(v)\,.
$$
This implies
$$
|D_\zs{1}(n)|\le \frac{1}{\gamma^{2}_\zs{j}}
\left(
3n |a|
+10
\right)\,.
$$
Similarly, one gets
$$
|D_\zs{2}(n)|\le \frac{2}{\gamma^{2}_\zs{j}}\,.
$$
Thus, for $j\ge 2$,
\begin{equation}\label{sec:Pr.2}
|\E_\zs{Q,S} \xi^2_\zs{j,n}
-
\varrho^{*}
 |
\le \frac{15\,|a|(1+|a|)\varrho^{*}}{\pi^{2}j^{2}}\,.
\end{equation}
Therefore
$$
\L_\zs{1,n}(Q)\le 2(1+|a|(|a|+1))\varrho^{*}
$$
and taking into account the conditions \eqref{sec:Ex.3} we get
\begin{equation}\label{sec:Pr.3}
\L_\zs{1,n}(Q)\le \L^{*}_\zs{1}\,,
\end{equation}
where $\L^{*}_\zs{1}$ is defined in \eqref{sec:Mo.15-1-1}. It means
that the condition $\C_\zs{1})$ holds with $\varsigma_\zs{Q}=\varrho^{*}$.
Moreover, by applying the conditions \eqref{sec:Ex.3} we have
$\varsigma^{*}=\varrho^{*}_\zs{max}$ and $\varsigma_\zs{*}=\varrho^{*}_\zs{min}$.

\noindent To check the condition $\C_\zs{2})$ we note that
$$
\E_\zs{Q,S}\, \left(
\sum^\infty_\zs{j=1}\,x_\zs{j}\,
(
\xi^2_\zs{j,n}-\E_\zs{Q,S} \xi^2_\zs{j,n}
)
\right)^{2}
=\frac{1}{n^{2}}\,\sum_\zs{i,j\ge 1} x_\zs{i} x_\zs{j}\,
\E_\zs{Q,S}\wt{I}_\zs{n}(\phi_\zs{i})\wt{I}_\zs{n}(\phi_\zs{j})\,.
$$
Therefore, in view of  Theorem~\ref{Th.sec:Ou.1}
\begin{equation}\label{sec:Pr.4}
\E_\zs{Q,S}\, \left(
\sum^\infty_\zs{j=1}\,x_\zs{j}\,
(
\xi^2_\zs{j,n}-\E_\zs{Q,S} \xi^2_\zs{j,n}
)
\right)^{2}
\le \frac{2\M^{*}}{n}
\sum_\zs{i,j\ge 1} |x_\zs{i}| |x_\zs{j}| (\varpi^{*}_\zs{i,j}+2)\,,
\end{equation}
where $\varpi^{*}_\zs{i,j}=\varpi^{*}_\zs{\phi_\zs{i},\phi_\zs{j}}$.
To estimate this term we note, that for any $j\ge 1$,
$$
\phi_\zs{j}(v+u)=a_\zs{j-1}(v)\phi_\zs{j-1}(u)
+a_\zs{j}(v)\phi_\zs{j}(u)
+a_\zs{j+1}(v)\phi_\zs{j+1}(u)
$$
$a_\zs{j}(\cdot)$ are bounded functions with $|a_\zs{j}(v)|\le 1$.
Thus,
$$
\varpi^{*}_\zs{i,j}\le 3n\Chi_\zs{\{|i-j|\le 1\}}+
3\Chi_\zs{\{|i-j|\ge 2\}}\,.
$$
Since
$\sum_\zs{j\ge 1} x^{2}_\zs{j}\le 1$,
therefore, the upper bound in \eqref{sec:Pr.4} can be estimated as
$$
\sum_\zs{i,j\ge 1} |x_\zs{i}| |x_\zs{j}| (\varpi^{*}_\zs{i,j}+2)\le
14 n\,.
$$
>From here, it follows that
\begin{equation}\label{sec:Pr.5}
\L_\zs{2,n}(Q)\le 28 \M^{*}\,.
\end{equation}
\noindent Hence Theorem~\ref{Th.sec:2.2}.

\endproof

\section{Robust asymptotic efficiency}\label{sec:Ef}

In this Section we show that the model selection procedure \eqref{sec:Mo.14},
\eqref{sec:Ga.4},
\eqref{sec:Si.1}, \eqref{sec:Ga.3}
for estimating $S$
in the model \eqref{sec:In.1}
 is asymptotically efficient with respect to the robust risk \eqref{sec:In.6}.
 We assume that the unknown function $S$ in the model
 \eqref{sec:In.1} belongs to the Sobolev ball
\begin{equation}\label{sec:Ef.1}
W^{k}_\zs{r}=\{f\in \,\cC^{k}_\zs{per}[0,1]
\,,\,\sum_\zs{j=0}^k\,\|f^{(j)}\|^2\le r\}\,,
 \end{equation}
where $r>0\,,\ k\ge 1$ are some  parameters,
$\cC^{k}_\zs{per}[0,1]$ is the set of
 $k$ times continuously differentiable functions
$f\,:\,[0,1]\to\bbr$ such that $f^{(i)}(0)=f^{(i)}(1)$ for all
$0\le i \le k$. The functional class $W^{k}_\zs{r}$ can be written as
an ellipsoid in $l_\zs{2}$, i.e.
 \begin{equation}\label{sec:Ef.2}
W^{k}_\zs{r}=\{f\in\,\cC^{k}_\zs{per}[0,1]\,:\,
\sum_\zs{j=1}^{\infty}\,a_\zs{j}\,\theta^2_\zs{j}\,\le r\}
 \end{equation}
where $a_\zs{j}=\sum^k_{i=0}\left(2\pi [j/2]\right)^{2i}$.

\noindent We denote by $Q_\zs{0}$ the distribution of Winer process with the scale parameter
$\varsigma^{*}$ defined in \eqref{sec:Mo.15-3}.

\medskip

\noindent $\H_\zs{1})$
{\em Assume the distribution $Q_\zs{0}$ belongs to the family $\cQ_\zs{n}$.
}

\medskip

In this Section we will show that
the Pinsker constant for the robust risk \eqref{sec:In.6}
is given by the equation
\begin{equation}\label{sec:Ef.3}
R^{*}_\zs{k}=\left((2k+1)r\right)^{1/(2k+1)}\,
\left(
\frac{\varsigma^{*} k}{(k+1)\pi} \right)^{2k/(2k+1)}\,.
\end{equation}
It is well known that
 the optimal (minimax) rate
for the Sobolev ball $W^{k}_\zs{r}$
 is $n^{2k/(2k+1)}$ (see, for example, \cite{Pi}, \cite{Nu}).

We will see that asymptotically the robust risk
 \eqref{sec:In.6} normalized by this rate is bounded
from below by $R^{*}_\zs{k}$, i.e. this bound can not be diminished if one
considers the class of all admissible estimates for $S$.
Let $\Pi_\zs{n}$ be the set of all estimators $\wh{S}_\zs{n}$
measurable with respect to the sigma-algebra
$\sigma\{y_\zs{t}\,,\,0\le t\le n\}$
 generated by the process \eqref{sec:In.1}.

\begin{theorem}\label{Th.sec:Ef.1} Under the condition $\H_\zs{1})$
 \begin{equation}\label{sec:Ef.4}
\liminf_{n\to\infty}\,n^{2k/(2k+1)}\, \inf_{\wh{S}_\zs{n}\in\Pi_\zs{n}}\,\,
\sup_\zs{S\in W^{k}_\zs{r}} \,\cR^{*}_\zs{n}(\wh{S}_\zs{n},S) \ge
R^{*}_\zs{k} \,.
 \end{equation}
\end{theorem}
\noindent Proof of this theorem follows directly from
Theorem 3.2 in \cite{KoPe4}.

\noindent Now we show that,
under some conditions, the normalized robust risk for the model selection
procedure
is bounded from above by the same constant $R^{*}_\zs{k}$.

\begin{theorem}\label{Th.sec:Ef.2}
Assume that, in model \eqref{sec:In.1}, for each $n\ge 1$ the distribution
of $(\xi_\zs{t})_\zs{0\le t\le n}$ belongs to the family
$\cQ_\zs{n}$ satisfying the conditions
$\H_\zs{0})$.
 Then  the
robust risk \eqref{sec:In.6}
of the  model selection procedure $\wh{S}_\zs{*}$ defined in
\eqref{sec:Ga.4},
\eqref{sec:Si.1}, \eqref{sec:Ga.3}
 has
 the following asymptotic upper bound
 \begin{equation}\label{sec:Ef.5}
\limsup_\zs{n\to\infty}\,n^{2k/(2k+1)}\,
 \sup_\zs{S\in W^k_r}\,
\cR^{*}_\zs{n}(\wh{S}_\zs{*},S) \le  R^{*}_\zs{k}\,.
 \end{equation}
\end{theorem}
\medskip
\noindent Theorem~\ref{Th.sec:Ef.1} and Theorem~\ref{Th.sec:Ef.2} imply the following result
\begin{corollary}\label{Co.sec:Mr.1}
Under the conditions $\H_\zs{0})$ and
$\H_\zs{1})$
\begin{equation}\label{sec:Ef.6}
\lim_{n\to\infty}\,n^{2k/(2k+1)}\, \inf_{\wh{S}_\zs{n}\in\Pi_\zs{n}}\,\,
\sup_\zs{S\in W^{k}_\zs{r}} \,\cR^{*}_\zs{n}(\wh{S}_\zs{n},S)
= R^{*}_\zs{k}\,.
 \end{equation}
\end{corollary}

\begin{remark}\label{Re.sec:Mr.1}
The equation \eqref{sec:Ef.6} means that the parameter $R^{*}_\zs{k}$ defined by \eqref{sec:Ef.3}
is the Pinsker constant (see, for example, \cite{Pi}, \cite{Nu}) for the model \eqref{sec:In.1}.
Moreover, the equality \eqref{sec:Ef.6} means  that the
model selection procedure
\eqref{sec:Ga.4},
\eqref{sec:Si.1}, \eqref{sec:Ga.3}
is asymptotically robust efficient.
\end{remark}

\section{Upper bound}\label{sec:Up}

\subsection{Known smoothness}

First we suppose that the parameters $k\ge 1$, $r>0$ in \eqref{sec:Ef.1} and
$\varsigma^{*}$ in \eqref{sec:Mo.15-3}  are known. Let the family of admissible
weighted least squares estimates
$(\wh{S}_\zs{\gamma})_\zs{\gamma\in\Gamma}$
for the unknown function $S\in W^{k}_\zs{r}$ be given
by
\eqref{sec:Ga.3}. Consider the pair
$$
\alpha_\zs{0}=(k,t_\zs{0})
$$
where
 $t_\zs{0}=[\ov{r}/\varepsilon]\varepsilon$,
$\ov{r}=r/\varsigma^{*}$
and $\varepsilon$ satisfies the conditions
in \eqref{sec:Ga.1}.
Denote the corresponding weight sequence in $\Gamma$ as
\begin{equation}\label{sec:Up.1}
\gamma_\zs{0}=\gamma_\zs{\alpha_\zs{0}}\,.
\end{equation}
\noindent Note that for sufficiently large $n$ the pair
$\alpha_\zs{0}$ belongs to the set \eqref{sec:Ga.0}.

\begin{theorem}\label{Th.sec:Up.1}
The estimator $\wh{S}_\zs{\gamma_\zs{0}}$ satisfies the following asymptotic upper
bound
\begin{equation}\label{Sec:Up.2}
\limsup_{n\to\infty}\,n^{2k/(2k+1)}\,
\sup_{S\in W^{k}_\zs{r}}\,
\cR^{*}_\zs{n}\,(\wh{S}_\zs{\gamma_\zs{0}},S)\,
\le R^{*}_\zs{k}\,.
\end{equation}
\end{theorem}
\noindent {\bf Proof.}
Substituting the model \eqref{sec:In.1} in the definition of
$\wh{\theta}_\zs{j,n}$ in \eqref{sec:Mo.4}
yields
$$
\wh{\theta}_\zs{j,n}=
\theta_\zs{j}
+
\frac{1}{\sqrt{n}}\xi_\zs{j,n}\,,
$$
where the random variables $\xi_\zs{j,n}$ are defined in
\eqref{sec:Mo.4}. Therefore, by the definition of
the estimators $\wh{S}_\zs{\gamma}$ in \eqref{sec:Mo.6},
we get
$$
\|\wh{S}_\zs{\gamma_\zs{0}}-S\|^2
=\sum_{j=1}^{n}\,(1-\gamma_\zs{0}(j))^2\,\theta^2_\zs{j}-2M_\zs{n}
 +\,\frac{1}{n}\,
 \sum_{j=1}^{n}\,\gamma_\zs{0}^2(j)\,\xi^2_\zs{j,n}
$$
with
$$
M_\zs{n}\,=\,\frac{1}{\sqrt{n}}
\sum_{j=1}^{n}\,(1\,-\,\gamma_\zs{0}(j))\,\gamma_\zs{0}(j)\,\theta_\zs{j}\,\xi_\zs{j,n}
\,.
$$
It should be observed that
$\E_\zs{Q,S}\,M_\zs{n}=0$
for any $Q\in\cQ^{*}_\zs{n}$.  Moreover, by the condition $\C_\zs{1})$
$$
\E_\zs{Q,S}\,\sum_{j=1}^{n}\,\gamma_\zs{0}^2(j)\xi^{2}_\zs{j,n}
\le \varsigma_\zs{Q}\sum_{j=1}^{n}\,\gamma_\zs{0}^2(j)+ \L_\zs{1,n}(Q)
$$
and, taking into account the condition $\H_\zs{0})$, we get
$$
\sup_\zs{Q\in\cQ_\zs{n}}\,
\E_\zs{Q,S}
\sum_{j=1}^{n}\,\gamma_\zs{0}^2(j)
\,\xi^{2}_\zs{j,n}\le \varsigma^{*}\sum_{j=1}^{n}\,\gamma_\zs{0}^2(j)
+l_\zs{n}\,.
$$
Thus,
\begin{equation}\label{sec:Up.3}
\cR^{*}_\zs{n}(\wh{S}_\zs{\gamma_\zs{0}},S)
\,\le\,
\sum_{j=\iota_\zs{0}}^{n}\,(1-\gamma_\zs{0}(j))^2\,\theta^2_\zs{j}+
\frac{\varsigma^{*}}{n}\sum_{j=1}^{n}\,\gamma_\zs{0}^2(j)
+\frac{l_\zs{n}}{n}
\end{equation}
where $\iota_\zs{0}=j_\zs{0}(\alpha_\zs{0})$. Setting
$$
\upsilon_\zs{n}= n^{2k/(2k+1)} \sup_\zs{j\ge
\iota_\zs{0}}(1-\gamma_\zs{0}(j))^2/a_\zs{j}\,,
$$
we estimate the first summand in the right-hand of \eqref{sec:Up.3}
as
$$
n^{2k/(2k+1)}\,
\sum_{j=\iota_\zs{0}}^{n}\,
(1-\gamma_\zs{0}(j))^2\,\theta^2_\zs{j}
\le
\upsilon_\zs{n}\,
\sum_\zs{j\ge 1}\,a_\zs{j}\,\theta^{2}_\zs{j}
\,.
$$
>From here and
\eqref{sec:Ef.2}, we obtain that for each $S\in W^{k}_\zs{r}$
$$
\Upsilon_\zs{1,n}(S)=\,n^{2k/(2k+1)}\,
\sum_{j=\iota_\zs{0}}^{n}\,(1-\gamma_\zs{0}(j))^2\,\theta^2_\zs{j}
\,\le\,\upsilon_\zs{n}\,r\,.
$$
Further we note that
$$
\limsup_{n\to\infty}\,\left(\ov{r}\right)^{2k/(2k+1)}\,
\upsilon_\zs{n}
\le\,
\frac{1}{\pi^{2k}\left(\tau_\zs{k}\right)^{2k/(2k+1)}}
\,,
$$
where the coefficient $\tau_\zs{k}$ is given in \eqref{sec:Ga.2}.
Therefore,
\begin{equation}\label{Up.4}
\limsup_\zs{n\to\infty}
\sup_\zs{S\in W^{k}_\zs{r}}\,
\Upsilon_\zs{1,n}(S)
\,\le\,
(\varsigma^{*})^{2k/(2k+1)}\,
\Upsilon^{*}_\zs{1}
\end{equation}
where
$$
\Upsilon^{*}_\zs{1}=
\frac{r^{1/(2k+1)}}{\pi^{2k}(\tau_\zs{k})^{2k/(2k+1)}}\,.
$$
To examine the second summand in the right hand of \eqref{Sec:Up.2}, we set
$$
\Upsilon_\zs{2,n}=\frac{1}{n^{1/(2k+1)}}\sum_{j=1}^{n}\,\gamma_\zs{0}^2(j)\,.
$$
It is easy to check that
$$
\lim_{n\to\infty}\,
\frac{1}{(\ov{r})^{1/(2k+1)}}\,
\Upsilon_\zs{2,n}=
\frac{2(\tau_\zs{k})^{1/(2k+1)}\,k^2}{(k+1)(2k+1)}
:=\Upsilon^{*}_\zs{2}
\,.
$$
Therefore, taking into account that
$$
(\varsigma^{*})^{2k/(2k+1)}\,
\Upsilon^{*}_\zs{1,n}
+
\varsigma^{*}(\ov{r})^{1/(2k+1)}
\Upsilon^{*}_\zs{2}
=
R^{*}_\zs{k}\,,
$$
we obtain
$$
\lim_\zs{n\to\infty}
n^{2k/(2k+1)}\,\sup_\zs{S\in W^{k}_\zs{r}}
\,\cR^{*}_\zs{n}(\wh{S}_\zs{\gamma_\zs{0}},S)
\le\,R^{*}_\zs{k}\,.
$$
Hence Theorem~\ref{Th.sec:Up.1}.
\endproof

\subsection{Unknown smoothness}

Combining Theorem~\ref{Th.sec:Up.1} and Theorem~\ref{Th.sec:La.1}
yields Theorem~\ref{Th.sec:Ef.2}.

\endproof

\renewcommand{\theequation}{A.\arabic{equation}}
\renewcommand{\thetheorem}{A.\arabic{theorem}}
\renewcommand{\thesubsection}{A.\arabic{subsection}}
\section{Appendix}\label{sec:A}
\setcounter{equation}{0}
\setcounter{theorem}{0}

\subsection{Technical lemmas}\label{subsec:A.2}

\begin{lemma}\label{Le.sec:A.1}
The operators $\tau_\zs{f,g}$ and $H_\zs{f,g}$ satisfy the following inequalities
\begin{equation}
\label{sec:A.1}
\sup_\zs{0\le t\le n}\,|\tau_\zs{f,g}(t)|\,\le 4 \varpi^{*}_\zs{f,g}
\quad\mbox{and}
\quad
\sup_\zs{0\le t\le n}\,|H_\zs{f,g}(t)|\,\le \D^{*}_\zs{1} \varpi^{*}_\zs{f,g}
\,,
\end{equation}
where
$\D^{*}_\zs{1}$ is given in \eqref{sec:Mo.15-1-2}.
\end{lemma}
\proof
Fist note that
$$
\int^{t}_\zs{0} f(s)  \varepsilon_\zs{g}(s) \d s=
a\int^{t}_\zs{0}\,e^{a v}
\left(
\int^{t-v}_\zs{0}
f(s+v)g(s)(1+e^{2as})\d s
\right)\,\d v\,.
$$
Integrating by parts yields
\begin{align*}
\int^{t-v}_\zs{0}
f(s+v)g(s)(1+e^{2as})\d s
&=
(1+e^{2a(t-v)}) \int^{t-v}_\zs{0} f(s+v) g(s)\d s\\[2mm]
&-
2a \int^{t-v}_\zs{0} e^{2as}
\left(
 \int^{s}_\zs{0} f(z+v) g(z)\d z
\right)
\d s\,.
\end{align*}
Taking into account the definition \eqref{sec:Ou.10}, we estimate this integral as
$$
\left|
\int^{t-v}_\zs{0}
f(s+v)g(s)(1+e^{2as})\d s
\right|\le 3 \varpi^{*}_\zs{f,g}\,.
$$
Therefore,
$$
\left|
\int^{t}_\zs{0} f(s)  \varepsilon_\zs{g}(s) \d s
\right|
\le 3 \varpi^{*}_\zs{f,g}
\quad\mbox{and}\quad
\left|
\int^{t}_\zs{0} \varepsilon^{*}_\zs{f,g}(s) \d s
\right|
\le 6 \varpi^{*}_\zs{f,g}
\,.
$$
This implies the first inequality in \eqref{sec:A.1}.
To obtain the second one we represent the function $H_\zs{f,g}(t)$ in the following
form
\begin{align*}
H_\zs{g,f}(t)&=\lambda \varrho^{2}_\zs{1} \tau_\zs{f,g}(t)+
\lambda^{2} \varrho^{2}_\zs{2} \int^{t}_\zs{0} f(z) g(z) \d z\\[2mm]
&+
\lambda^{2} \varrho^{2}_\zs{2}
\left(H^{(1)}_\zs{g,f}(t)+H^{(1)}_\zs{f,g}(t)
+H^{(2)}_\zs{g,f}(t)+H^{(2)}_\zs{f,g}(t)
\right)\,,
\end{align*}
where
$$
H^{(1)}_\zs{g,f}(t)=a\int^{t}_\zs{0}\,\int^{t-z}_\zs{0}\,e^{ay}
g(y+z)f(z)\,\d y\,\d z
$$
and
$$
H^{(2)}_\zs{g,f}(t)=a^{2}\int^{t}_\zs{0}\,\int^{t-z}_\zs{0}\,e^{ay}
g(y+z)
\int^{y}_\zs{0} e^{av} f(v+z)\,
\d v\,\d y\,\d z\,.
$$
Now we note
$$
|H^{(1)}_\zs{g,f}(t)|=
\left|
a\int^{t}_\zs{0}\,e^{ay}\,
\left(
\int^{t-y}_\zs{0}\,
g(y+z)f(z)\,\d z \right)\d y
\right|\le \varpi^{*}_\zs{f,g}\,.
$$
To estimate $H^{(2)}_\zs{g,f}(t)$ we represent it as
$$
H^{(2)}_\zs{g,f}(t)=a^{2}\int^{t}_\zs{0} e^{a y}\,
\left(
\int^{y}_\zs{0} e^{av}
\left(
\int^{t-y}_\zs{0}
g(y+z)f(v+z)\d z
\right)
\d v
\right)
\d y\,.
$$
Note that
for any $0\le v\le y\le t$
one has
$$
\left|
\int^{t-y}_\zs{0}
g(y+z)f(v+z)\d z
\right|
%=
%\left|
%\int^{t-y+v}_\zs{v}
%g(z+y-v)f(z)\d z
%\right|\\[2mm]
%&=
%\left|
%\int^{t-y+v}_\zs{0}
%g(z+y-v)f(z)\d z
%-
%\int^{v}_\zs{0}
%g(z+y-v)f(z)\d z
%\right|
\le 2\varpi^{*}_\zs{f,g}\,.
$$
Thus, $|H^{(2)}_\zs{g,f}(t)|\le 2\varpi^{*}_\zs{f,g}$,
and we come to the second inequality in \eqref{sec:A.1}.
Hence Lemma~\ref{Le.sec:A.1}.

\endproof

\medskip

\begin{lemma}\label{Le.sec:A.2}
For any bounded left-continuous
$[0,+\infty)\to \bbr$ functions $f$, $g$
$$
\E\,v_\zs{t}(f)\,v_\zs{t}(g)=
\int^{t}_\zs{0}\,e^{2a(t-s)}\,
V_\zs{f,g}(s)
\,
\d s\,,
$$
where $V_\zs{f,g}(s)$ is given in \eqref{sec:Ou.19}.
\end{lemma}
\proof
By the Ito formula and \eqref{sec:Ou.14}, one gets
\begin{align*}
\d \E v_\zs{t}(f) v_\zs{t}(g)&=2a \E v_\zs{t}(f) v_\zs{t}(g)\d t
+
(G^{*}_\zs{f,g}(t)+\varrho_\zs{3}\,f(t)g(t))
\d t\\[2mm]
&+a\left(
g(t)
\E v_\zs{t}(f)\,\zeta_\zs{t}
+f(t)
\E v_\zs{t}(g)\,\zeta_\zs{t}
\right)\d t\,.
\end{align*}
To calculate $\E v_\zs{t}(f)\,\zeta_\zs{t} $, we put $g=1$ in this equality.
Then, taking into account that
$$
\kappa_\zs{f}(t)=aG^{*}_\zs{f,1}(t)+a\varrho_\zs{3}\,f(t)\,,
$$
we get
\begin{equation}\label{sec:A.2}
a\E v_\zs{t}(f)\,\zeta_\zs{t} =
\int^{t}_\zs{0}\,e^{3a(t-s)}
\left(f(s) a^{2} \E\,\zeta^{2}_\zs{s}
+\kappa_\zs{f}(s)
\right)\d s
=\,A_\zs{f}(t)\,.
\end{equation}
Therefore
\begin{align*}
\E v_\zs{t}(f) v_\zs{t}(g)&=
\int^{t}_\zs{0}\,e^{2a(t-s)}
\left(
g(s)\,A_\zs{f}(s)
+f(s)
\,A_\zs{g}(s)
\right)\d s\\[2mm]
&+\int^{t}_\zs{0}\,e^{2a(t-s)}
\left(
G^{*}_\zs{f,g}(s)+\varrho_\zs{3}\,f(s)g(s)
\right)
\d s\,.
\end{align*}
Hence Lemma~\ref{Le.sec:A.2}.

\endproof

\medskip

\medskip
\noindent Further we will need the following result.
\begin{lemma}\label{Le.sec:A.3}
Let $\upsilon$ be a continuously differentiable $\bbr\to\bbr$ function. Then, for any
$n\ge 1$, $\alpha>0$ and for any  integrated $\bbr\to\bbr$ function $\Psi$,
$$
\sup_\zs{0\le t\le n}\,
\left|\int^{t}_\zs{0}\,e^{-\alpha(t-s)} \Psi(s) \upsilon(s)\,\d s \right|\,
\le \,
\varpi_\zs{1,\Psi}
\left(
2 \|\upsilon\|_\zs{*}
+\frac{\|\dot{\upsilon}\|_\zs{*}}{\alpha}\,
\right)\,.
$$
\end{lemma}
\proof
One obtains this inequality with the help of integrating by parts.

\endproof

\medskip

\begin{lemma}\label{Le.sec:A.4}
For any mesurable bounded  $[0,+\infty)\to\bbr$ functions $f$ and $g$,
for any $-\infty< a\le 0 $ , for any $ 0\le t\le n$ and for any $n\ge 1$
\begin{equation}\label{sec:A.3}
\left|
a
\int^{t}_\zs{0}
e^{2a(t-s)}\,g(s)
A_\zs{f}(s)
\,\d s
\right|\,\le
3\lambda_\zs{2}
\varpi^{*}_\zs{f,g}
+
\varrho_\zs{3}\|f\|_\zs{*}
\|g\|_\zs{*}
\,.
\end{equation}
\end{lemma}
\proof
One can represent the function $A_\zs{f}(t)$ as
\begin{equation}\label{sec:A.4}
A_\zs{f}(t)=
\int^{t}_\zs{0} e^{3a(t-s)}f(s)\upsilon(s)\d s
+
\lambda_\zs{2}
\int^{t}_\zs{0} e^{3a(t-s)}\varepsilon_\zs{f}(s)\d s
\,,
\end{equation}
where $\upsilon(s)=a^{2}\E\,\zeta^{2}_\zs{s}+\lambda_\zs{2}
\left(
e^{2as}-1
\right)
+a\varrho_\zs{3}$.
>From here and  \eqref{sec:Ou.22} we rewrite this function as
\begin{equation}\label{sec:A.5}
\upsilon(s)=
a\upsilon_\zs{1}(s)
+
\upsilon_\zs{2}(s)
\end{equation}
with
$$
\upsilon_\zs{1}(s)=
\frac{\varrho_\zs{3}}{4}
\left(
e^{4at}+3
\right)
\quad\mbox{and}\quad
\upsilon_\zs{2}(s)=
\frac{\lambda_\zs{2}}{2}\,
e^{4at}
-
\frac{\lambda_\zs{2}}{2}
\,.
$$
These functions can be estimated as
\begin{equation}\label{sec:A.6}
\begin{array}{rl}
\|\upsilon_\zs{1}\|_\zs{*}
&\le\,\varrho_\zs{3}\,;\\[4mm]
\sup_\zs{-\infty< a\le 0}\,
\left(
2\|\upsilon_\zs{2}\|_\zs{*}
+
\dfrac{\|\dot{\upsilon}_\zs{2}\|_\zs{*}}{2|a|}
\right)\,
&\le\,2\lambda_\zs{2}\,.
\end{array}
\end{equation}
Now we represent the intergal in \eqref{sec:A.3} as
$$
a
\int^{t}_\zs{0}
e^{2a(t-s)}\,g(s)
A_\zs{f}(s)
\,\d s
=J_\zs{1}(t)
+
J_\zs{2}(t)
+\lambda_\zs{2}J_\zs{3}(t)\,,
$$
where
\begin{align*}
 J_\zs{1}(t)&=
a^{2}
\int^{t}_\zs{0}
e^{2a(t-s)}\,g(s)
\left(
\int^{s}_\zs{0}e^{3a(s-u)} f(u) \upsilon_\zs{1}(u)
\d u \right)\,
\d s\,,\\[2mm]
J_\zs{2}(t)&=
a
\int^{t}_\zs{0}
e^{2a(t-s)}\,g(s)
\left(
\int^{s}_\zs{0}e^{3a(s-u)} f(u) \upsilon_\zs{2}(u)
\d u \right)\,
\d s\,,\\[2mm]
J_\zs{3}(t)&=a
\int^{t}_\zs{0}
e^{2a(t-s)}\,g(s)
\left(
\int^{s}_\zs{0}e^{3a(s-u)} \varepsilon_\zs{f}(u)
\d u \right)\,
\d s\,.
\end{align*}
In view of \eqref{sec:A.6} we have
$$
\sup_\zs{0\le t\le n}
|J_\zs{1}(t)|\le \varrho_\zs{3}
\|f\|_\zs{*} \|g\|_\zs{*}
\,.
$$
Further we represent $J_\zs{2}(t)$ as
$$
J_\zs{2}(t)=a\int^{t}_\zs{0}\,e^{3a u}\,
\left(
\int^{t-u}_\zs{0}\,e^{2a(t-u-s)}\,g(s+u) f(s)\upsilon_\zs{1}(s) \d s
\right)
 \d u\,.
$$
By Lemma~\ref{Le.sec:A.3} and \eqref{sec:A.6}
 we obtain that for any $0\le z\le n$ and $0\le u\le n-z$
$$
\left|\int^{z}_\zs{0} e^{2a(z-s)} \upsilon(s) g(s+u) f(s)\d s\right|\le
3\lambda_\zs{2}\varpi^{*}_\zs{f,g} \,.
$$
Therefore,
$$
\sup_\zs{0\le t\le n}
|J_\zs{2}(t)|
\le \lambda_\zs{2}\,\varpi^{*}_\zs{f,g}\,.
$$
Similarly, one gets
$$
\sup_\zs{0\le t\le n}
|J_\zs{3}(t)|
\le \frac{2}{3}\,\varpi^{*}_\zs{g,\varepsilon_\zs{f}}\,.
$$
To estimate the quantity $\varpi_\zs{g,\varepsilon_\zs{f}}$
defined in \eqref{sec:Ou.10} we
 note that
 for any   $0\le v\le n$ and $0\le t\le n-v$
\begin{equation}\label{sec:A.6'}
\int^{t}_\zs{0}\,g(s+v)\varepsilon_\zs{f}(s)\d s=
a\int^{t}_\zs{0}\,e^{ax}\,\Theta_\zs{g,f}(t-x,v+x)\d x\,,
\end{equation}
where
$$
\Theta_\zs{g,f}(t,v)=\int^{t}_\zs{0}\,g(s+v)f(s)(1+e^{2as})\d s\,.
$$
Denoting
\begin{equation}\label{sec:A.7}
\Upsilon_\zs{g,f}(s,u)=\int^{s}_\zs{0}\,g(r+u) f(r)\d r\,,
\end{equation}
we represent the function $\Theta_\zs{g,f}(t,v)$ as
$$
\Theta_\zs{g,f}(t,v)=(1+e^{2at})\,\Upsilon_\zs{g,f}(t,v)
-2a\int^{t}_\zs{0}\,e^{2as}\Upsilon_\zs{g,f}(s,v)\d s\,.
$$
Therefore
$$
\max_\zs{0\le v\le n}\,
\max_\zs{0\le t\le n-v}\,
|
\Theta_\zs{g,f}(t,v)
|\,\le \,3\varpi^{*}_\zs{f,g}\,.
$$
In view of \eqref{sec:A.6'}, one gets
$$
\varpi_\zs{g,\varepsilon_\zs{f}}
\le 3\varpi^{*}_\zs{f,g}
\quad\mbox{and}\quad
\sup_\zs{0\le t\le n}
|J_\zs{3}(t)| \le 2 \varpi^{*}_\zs{f,g}\,.
$$
Hence Lemma~\ref{Le.sec:A.4}.

\endproof

\medskip

\begin{lemma}\label{Le.sec:A.5}
For any mesurable bounded  $[0,+\infty)\to\bbr$ functions $f$ and $g$,
for any $-\infty< a\le 0 $ , for any $ 0\le t\le n$ and for any $n\ge 1$
\begin{equation}\label{sec:A.8}
\left|
2a
\int^{t}_\zs{0}
e^{2a(t-s)}\,
G^{*}_\zs{f,g}(s)
\,\d s
\right|\,\le
\,\D^{*}_\zs{2}\,\varpi_\zs{f,g}
\,,
\end{equation}
where
$\D^{*}_\zs{2}$ is defined in
\eqref{sec:Mo.15-1-2}.
\end{lemma}
\proof
First we note that
the function $G^{*}_\zs{f,g}$ can be represented  as
$$
G^{*}_\zs{f,g}(t)=G_\zs{f,g}(t)+\frac{\lambda_\zs{2}}{a}\,
\dot{\tau}_\zs{f,g}(t)
+
\frac{\lambda_\zs{2}}{2a}
\left(
e^{2at}-3
\right)
f(t)g(t)\,.
$$
Integrating by the parts yields
$$
\left|
\int^{t}_\zs{0}
e^{2a(t-s)}\,
\dot{\tau}_\zs{f,g}(s)
\,\d s
\right|\,\le 8 \varpi^{*}_\zs{f,g}\,.
$$
Finally, by applying Lemma~\ref{Le.sec:A.3} with $\upsilon(s)=e^{2as}-3$
and $\Psi(s)=f(s)g(s)$, one gets
$$
\left|
\int^{t}_\zs{0}
e^{2a(t-s)}\,
(e^{2as}-3)f(s)g(s)
\,\d s
\right|\,\le 7 \varpi^{*}_\zs{f,g}\,.
$$
\noindent Hence Lemma~\ref{Le.sec:A.5}.

\endproof

\begin{lemma}\label{Le.sec:A.6}
For any mesurable bounded  $[0,+\infty)\to\bbr$ functions $f$ and $g$,
for any $-\infty< a\le 0 $ , for any $ 0\le t\le n$ and for any $n\ge 1$
\begin{equation}\label{sec:A.9}
|\cV_\zs{f,g}(t)|
\le\,
(6\lambda_\zs{2}+\D^{*}_\zs{2})\varpi^{*}_\zs{f,g}
+3\varrho_\zs{3}\|f\|_\zs{*} \|g\|_\zs{*}
\,.
\end{equation}
\end{lemma}
\proof
This inequality is a direct consequence of Lemmas~\ref{Le.sec:A.2}-\ref{Le.sec:A.5}
Hence Lemma~\ref{Le.sec:A.6}.

\endproof

\medskip

\begin{lemma}\label{Le.sec:A.7}
For any mesurable bounded  $[0,+\infty)\to\bbr$ functions $f$ and $g$,
for any $-\infty< a\le 0 $ , $ 0\le t\le n$ and $n\ge 1$
\begin{equation}\label{sec:A.10}
\left|
a
\int^{t}_\zs{0}
e^{a(t-s)}\,f(s)
\cK_\zs{g}(s)
\,\d s
\right|\le 14 \lambda_\zs{2}
\|g\|_\zs{*}
\varpi^{*}_\zs{f,g}+\varrho_\zs{3}\|f\|_\zs{*}\,\|g\|^{2}_\zs{*}\,.
\end{equation}
\end{lemma}
\proof
Taking into account \eqref{sec:A.4}--\eqref{sec:A.5}, we write down the function
$A_\zs{g}(t)$ as
$$
A_\zs{g}(t)=A^{(1)}_\zs{g}(t)+A^{(2)}_\zs{g}(t)
$$
where
$$
A^{(1)}_\zs{g}(t)=a\int^{t}_\zs{0} e^{3a(t-s)}g(s)\upsilon_\zs{1}(s)\d s
\,, \quad
A^{(2)}_\zs{g}(t)=
\int^{t}_\zs{0} e^{3a(t-s)}
\left(
g(s)\upsilon_\zs{2}(s)
+
\lambda_\zs{2}
\varepsilon_\zs{g}(s)
\right)
\d s
\,.
$$
Since
$$
U_\zs{1,g}(t)=\frac{2\lambda_\zs{2}}{a}g(t)\varepsilon_\zs{g}(t)+
\varrho_\zs{3}g^{2}(t)\,,
$$
 the integral in \eqref{sec:A.11} can be represented as
$$
a
\int^{t}_\zs{0}
e^{a(t-s)}\,f(s)
\cK_\zs{g}(s)
\,\d s
=J^{*}_\zs{1}(t)+J^{*}_\zs{2}(t)+J^{*}_\zs{3}(t)
+J^{*}_\zs{4}(t)
\,,
$$
where
\begin{align*}
J^{*}_\zs{1}(t)&=2a^{2}\int^{t}_\zs{0}e^{a(t-s)}f(s)\,
\left(\int^{s}_\zs{0} e^{2a(s-r)} g(r) A^{(1)}_\zs{g}(r)\d r
\right)\d s\,,\\[2mm]
J^{*}_\zs{2}(t)&=2a^{2}\int^{t}_\zs{0}e^{a(t-s)}f(s)\,
\left(\int^{s}_\zs{0} e^{2a(s-r)} g(r) A^{(2)}_\zs{g}(r)\d r
\right)\d s\,,\\[2mm]
J^{*}_\zs{3}(t)&=2a\lambda_\zs{2}\int^{t}_\zs{0}e^{a(t-s)}f(s)\,
\left(
\int^{s}_\zs{0} e^{2a(s-r)} g(r) \varepsilon_\zs{g}(r)\d r
\right)\d s\,,\\[2mm]
J^{*}_\zs{4}(t)&=\varrho_\zs{3}\,a^{2}\int^{t}_\zs{0}e^{a(t-s)}f(s)\,\left(
\int^{s}_\zs{0} e^{2a(s-r)} g^{2}(r) \d r\right)\d s\,.
\end{align*}
In view of \eqref{sec:A.6}, one obtains
$$
\sup_\zs{0\le t\le n}
|A^{(1)}_\zs{g}(t)|\le \frac{\varrho_\zs{3}}{3}\|g\|_\zs{*}
\quad\mbox{and}\quad
\sup_\zs{0\le t\le n}\,|J^{*}_\zs{1}(t)|\le \frac{\varrho_\zs{3}}{3} \|f\|_\zs{*}\,\|g\|^{2}_\zs{*}\,.
$$
Denoting
$$
\Gamma_\zs{f,g}(t,x)=\int^{t}_\zs{0} e^{a(t-s)} A^{(2)}_\zs{g}(s) f(s+x)g(s)\d s\,,
$$
one has
$$
J^{*}_\zs{2}(t)=2a^{2} \int^{t-x}_\zs{0} e^{2ax} \Gamma_\zs{f,g}(t-x,x)\,\d x\,.
$$
Noting that
$$
\sup_\zs{0\le t\le n}
|\upsilon_\zs{2}(t)|\le \frac{\lambda_\zs{2}}{2}\,,
$$
one comes to the inequalities
$$
\sup_\zs{0\le t\le n}
|A^{(2)}_\zs{g}(t)|\le
\,
\frac{5\lambda_\zs{2}}{6|a|}
\|g\|_\zs{*}
\quad\mbox{and}\quad
\sup_\zs{0\le t\le n}
|\dot{A}^{(2)}_\zs{g}(t)|\le \,4 \lambda_\zs{2}\,\|g\|_\zs{*}\,.
$$
By applying Lemma~\ref{Le.sec:A.3} with $\Psi(s)=f(s+x)g(s)$ and $\upsilon(s)=\upsilon_\zs{2}(s)$
one gets
$$
\sup_\zs{0\le t\le n}
\sup_\zs{0\le x\le t}
|\Gamma_\zs{f,g}(t-x,x)|\le \frac{17}{3|a|}\,\lambda_\zs{2}\,
\varpi^{*}_\zs{f,g}\,\|g\|_\zs{*}
\le
6 \lambda_\zs{2}\,
\varpi^{*}_\zs{f,g}\,\|g\|_\zs{*}\,.
$$
Therefore,
$$
\sup_\zs{0\le t\le n}
|J^{*}_\zs{2}(t)|\le
6 \lambda_\zs{2}\,
\varpi^{*}_\zs{f,g}\,\|g\|_\zs{*}
\,.
$$
Similarly, one can show that
$$
\sup_\zs{0\le t\le n}
|J^{*}_\zs{3}(t)|\le 8\,\lambda_\zs{2}\,\varpi^{*}_\zs{f,g}\,\|g\|_\zs{*}\,.
$$
Finally, the function $J^{*}_\zs{4}(t)$ can be estimated as
$$
\sup_\zs{0\le t\le n}
|J^{*}_\zs{4}(t)|\le \frac{\varrho_\zs{3}}{2}
\|f\|_\zs{*}\,\|g\|^{2}_\zs{*}\,.
$$
Hence Lemma~\ref{Le.sec:A.7}.

\endproof

\medskip

\medskip

\begin{lemma}\label{Le.sec:A.8}
For any measurable bounded  $[0,+\infty)\to\bbr$ functions $f$ and $g$,
$-\infty< a\le 0$,  $0\le t\le n$ and $n\ge 1$
\begin{align}\label{sec:A.11}
\left|
a
\int^{t}_\zs{0}
e^{a(t-s)}\,
U_\zs{f,g}(s)
\,\d s
\right| \le 2
\D^{*}_\zs{2}
\,\|g\|_\zs{*} \varpi^{*}_\zs{f,g}
+\varrho_\zs{3} \|f\|_\zs{*}\,\|g\|^{2}_\zs{*}
\,.
\end{align}
\end{lemma}
\proof
We note that
$$
U_\zs{f,g}(t)=
2g(t) G_\zs{f,g}(t)+
\frac{\lambda_\zs{2}}{a} f(t)g(t)\varepsilon_\zs{g}(t)+
\varrho_\zs{3}f(t)g^{2}(t)\,.
$$
Taking into account \eqref{sec:A.9} we obtain that
$$
\left|
2a
\int^{t}_\zs{0}
e^{a(t-s)}\,
g(s) G_\zs{f,g}(s)
\,\d s
\right|
\le
2(4\varrho^{2}_\zs{1}\varrho^{*}+\varrho_\zs{2}\D^{*}_\zs{1})\|g\|_\zs{*}\varpi^{*}_\zs{f,g}\,.
$$
In view of Lemma~\ref{Le.sec:A.3}, we obtain
$$
\sup_\zs{0\le t\le n}\,
\left|
\int^{t}_\zs{0}
e^{a(t-s)}\,
f(s)g(s) \varepsilon_\zs{g}(s)
\,\d s
\right|
\le \varpi^{*}_\zs{f,g}
\left(
2\|\varepsilon_\zs{g}\|_\zs{*}
+
\frac{\|\dot{\varepsilon}_\zs{g}\|_\zs{*}}{|a|}
\right)\,.
$$
Taking into account that
$$
\|\varepsilon_\zs{g}\|_\zs{*}\le 2\|g\|_\zs{*}
\quad\mbox{and}\quad
\|\dot{\varepsilon}_\zs{g}\|_\zs{*}\le 4|a|\|g\|_\zs{*}\,,
$$
one gets
$$
\sup_\zs{0\le t\le n}\,
\left|
\int^{t}_\zs{0}
e^{a(t-s)}\,
f(s)g(s) \varepsilon_\zs{g}(s)
\,\d s
\right|
\le 8\|g\|_\zs{*}\varpi^{*}_\zs{f,g}\,.
$$
>From here we come to desired result.
Hence lemma~\ref{Le.sec:A.8}.

\endproof

\medskip

\medskip

\subsection{Property of the Fourier coefficients}\label{subsec:A.3}
\begin{lemma}\label{Le.sec:A.9}
Suppose that the  function $S$ in \eqref{sec:In.1} is
differentiable and satisfies the condition \eqref{sec:Si.2}. Then
the Fourier coefficients \eqref{sec:Mo.2} satisfy the inequality
$$
\sup_\zs{l\ge 2}\, l\,\sum^\infty_\zs{j=l}\theta^2_\zs{j}\, \le\,
4\,|\dot{S}|^2_\zs{1}\,.
$$
\end{lemma}
{\bf Proof}. In view of \eqref{sec:Mo.1}, one has
$$
\theta_\zs{2p}=-
\frac{1}{\sqrt{2}\pi p}\int^1_\zs{0}
\,\dot{S}(t)\sin(2\pi pt)\d t
$$
and
\begin{align*}
\theta_\zs{2p+1}&=
\frac{1}{\sqrt{2}\pi p}\int^1_\zs{0}
\,\dot{S}(t)(\cos(2\pi pt)-1)\d t\\
&=
-
\frac{\sqrt{2}}{\pi p}\int^1_\zs{0}
\,\dot{S}(t)\sin^2(\pi pt)\d t\,,
\quad p\ge 1\,.
\end{align*}
>From here, it follows that
for any $j\ge 2$
$$
\theta^2_\zs{j}\,
\le\,
 \frac{2}{j^2}|\dot{S}|^2_\zs{1}\,.
$$
Taking into account that
$$
\sup_\zs{l\ge 2}
l\sum_\zs{j\ge l}\frac{1}{j^{2}}\,\le \,2\,,
$$
we arrive at the desired result.
\endproof

\section{Acknowledgments}
The paper is  supported by the RFBI-Grant 09-01-00172-a and
the Russian State Contract 02.740.11.5026.

\medskip

\medskip

%\begin{thebibliography}{100}

\end{document}